\documentclass[a4paper]{article}
\usepackage{amsmath,amssymb,amsthm,mathtools,mathrsfs}
\usepackage{graphicx}
\usepackage{url}
\usepackage{pdfsync}

\usepackage[utf8]{inputenc}
\usepackage[T1]{fontenc}
\usepackage[english]{babel}
\usepackage{lmodern}
\IfFileExists{stmaryrd.sty}{\usepackage{stmaryrd}\SetSymbolFont{stmry}{bold}{U}{stmry}{m}{n}}{%
\newcommand{\llbracket}{\mathopen{[\mkern-3mu[}}%
\newcommand{\rrbracket}{\mathclose{]\mkern-3mu]}}}
\usepackage{microtype}
\usepackage[margin=2.5cm]{geometry}
\usepackage{enumitem}
\usepackage[numbers,sort&compress]{natbib}
\usepackage{xcolor}
\usepackage{tikz}
\usetikzlibrary{arrows.meta}
\usepackage{needspace}
\usepackage{float}
\usepackage{hyperref}

\renewcommand{\geq}{\geqslant}
\renewcommand{\leq}{\leqslant}

\hypersetup{colorlinks=true,linkcolor=blue!60!black,citecolor=blue!60!black,urlcolor=blue!60!black,pdftitle={Lie-Lions brackets of nonlocal vector fields},pdfauthor={Vincent Boulard, Amaury Hayat, Benedetto Piccoli, Emmanuel Tr\'elat}}
\newcommand{\R}{\mathbb R}

\newcommand{\Ptwo}{\mathcal P_2(\mathbb R^d)}

\newcommand{\Id}{\mathrm{Id}}

\newcommand{\br}[2]{\llbracket #1,#2\rrbracket}
\newcommand{\SL}{\mathcal{S}_L^1}
\newcommand{\SLm}[1]{\mathcal{S}_{L,\leq #1}^1}

\newtheorem{theorem}{Theorem}[section]
\newtheorem{proposition}[theorem]{Proposition}
\newtheorem{lemma}[theorem]{Lemma}
\newtheorem{corollary}[theorem]{Corollary}
\theoremstyle{definition}
\newtheorem{definition}[theorem]{Definition}
\newtheorem{assumption}[theorem]{Assumption}
\theoremstyle{remark}
\newtheorem{remark}[theorem]{Remark}
\numberwithin{equation}{section}

\title{Lie-Lions brackets of nonlocal vector fields}

\author{
Vincent Boulard\thanks{CERMICS, CNRS, \'Ecole nationale des ponts et chauss\'ees, Institut Polytechnique de Paris, Marne-la-Vall\'ee, France (\texttt{vincent.boulard@enpc.fr}).}
\and
Amaury Hayat\thanks{CERMICS, CNRS, \'Ecole nationale des ponts et chauss\'ees, Institut Polytechnique de Paris, Marne-la-Vall\'ee, France (\texttt{amaury.hayat@enpc.fr}).}
\and
Benedetto Piccoli\thanks{Department of Mathematical Sciences and Center for Computational and Integrative Biology, Rutgers University, Camden, NJ, USA (\texttt{piccoli@camden.rutgers.edu}).}
\and
Emmanuel Tr\'elat\thanks{Sorbonne Universit\'e, Universit\'e Paris Cit\'e, CNRS, Inria, Laboratoire Jacques-Louis Lions, LJLL, F-75005 Paris, France (\texttt{emmanuel.trelat@sorbonne-universite.fr}).}
}

\date{}

\begin{document}

\maketitle

\begin{abstract}
We study nonlocal vector fields on the space of probability measures with finite second moment. Using Lions derivatives, we define their Lie-Lions bracket, which adds to the spatial Lie bracket a correction due to the evolution of the law, and prove that its vanishing is equivalent to commutation of the global flows, under first-order regularity and linear growth assumptions. We then construct a Lie algebra of complete polynomial nonlocal fields. For this class, we establish commutator and second-order Baker-Campbell-Hausdorff formulas with cubic remainders in the Wasserstein distance $W_2$, controlled by the second moment and uniform in the number of particles. The proofs rely on finite-particle realizations, which preserve the bracket and jointly determine the nonlocal field. We also obtain uniform estimates for smooth bounded kernels and a Leibniz rule for the action of polynomial fields on polynomial observables. Examples from collective dynamics show that the vanishing of the spatial bracket at a fixed law is neither necessary nor sufficient for commutation. They concern symmetries, aggregation-dependent rotation, synchronization, the transmission of a turn along a hierarchy, effective three-body interactions and occupancy-dependent circulation.
\end{abstract}



\section{Introduction}\label{sec:introduction}

Let $\Ptwo$ be the space of probability measures on $\R^d$ with finite second moment, endowed with the Wasserstein distance $W_2$ of order $2$. Given a nonlocal vector field $f:\R^d\times\Ptwo\to\R^d$, we consider the continuity equation
\begin{equation}\label{intro-continuity}
\partial_t\mu_t+\operatorname{div}\big(f(\cdot,\mu_t)\mu_t\big)=0.
\end{equation}
Such equations describe collective dynamics in which the velocity of an individual depends on the distribution of the population. They occur in models of aggregation \cite{CarrilloDiFrancescoFigalliLaurentSlepcev2011}, alignment \cite{CarrilloFornasierRosadoToscani2010,HaTadmor2008,MotschTadmor2011} and pedestrian motion \cite{BellomoPiccoliTosin2012,CarrilloMartinWolfram2016,ColGarLecPedestrian}. The measure formulation includes empirical distributions as well as continuous densities (see \cite{PicRosNonlocalVelocity}). We denote by $\Phi_t^f\mu$ the solution of \eqref{intro-continuity} at time $t$ with initial law $\mu$.

Our objective is to characterize commutation of two such evolutions and to estimate the effect of their order of composition. Since each evolution changes the law on which the other velocity depends, the spatial Lie bracket computed at a fixed law does not suffice. Informally, for regular fields $f$ and $g$, we have
$$
\Phi_{-t}^g\circ\Phi_{-t}^f\circ\Phi_t^g\circ\Phi_t^f(\mu)=\big(\Id+t^2\br{f}{g}(\cdot,\mu)\big)_\#\mu+\mathrm{O}(\vert t\vert^3)
$$
in $W_2$, where the nonlocal vector field $\br{f}{g}$, which we call the \emph{Lie-Lions bracket} of $f$ and $g$, is given by
$$
\br{f}{g}(x,\mu)=\underbrace{[f,g]_x(x,\mu)}_{\text{classical Lie bracket}}+\underbrace{\int_{\R^d}\big(\partial_\mu g(x,\mu)(y)f(y,\mu)-\partial_\mu f(x,\mu)(y)g(y,\mu)\big)\,d\mu(y)}_{\text{Lions correction }[f,g]_\mu(x,\mu)}.
$$
Here $[f,g]_x(x,\mu)=[f(\cdot,\mu),g(\cdot,\mu)](x)$ is the Lie bracket in $x$ at a fixed law, with the convention $[u,v]=Dv\,u-Du\,v$ for vector fields $u,v$ on a finite-dimensional space, and $\partial_\mu$ denotes the Lions derivative with respect to the law, recalled in Section~\ref{subsec:flows}. The matrix $\partial_\mu f(x,\mu)(y)$ measures how the velocity at $x$ varies when mass located near $y$ is moved. In both brackets, the index is the variable of differentiation. The Lions correction accounts for the change of the law under the two evolutions. Note that $[\cdot,\cdot]_\mu$ is bilinear and antisymmetric, but does not satisfy the Jacobi identity.\footnote{For $d=1$, $f(x,\mu)=x$, $g(x,\mu)=\int y\,d\mu(y)$ and $h=1$, we have $[[f,g]_\mu,h]_\mu+[[g,h]_\mu,f]_\mu+[[h,f]_\mu,g]_\mu=-1$.}

\paragraph{Commutation criterion.}
Under first-order regularity and linear growth assumptions (Assumption~\ref{ass:regularity}), the flows of $f$ and $g$ commute for all times and all laws if and only if $\br{f}{g}=0$ (Theorem~\ref{thm:commutation}). The proof uses the particle realizations $f_i^N(\boldsymbol{x})=f(x_i,\mu_{\boldsymbol{x}}^N)$, where $\mu_{\boldsymbol{x}}^N=\frac{1}{N}\sum_{i=1}^N\delta_{x_i}$. We prove that
$$
[f^N,g^N]=\br{f}{g}^N\qquad\forall N\geq1
$$
(Proposition~\ref{prop:particle}). This identity is exact for every $N$: the Lions correction is therefore not added to the mean-field model, but imposed by the bracket of the particle systems. Note that $[f^N,g^N]$ differs from $([f,g]_x)^N$ in general, since each component of $f^N$ depends on all positions through the empirical law. These realizations jointly determine a continuous nonlocal field (Lemma~\ref{lem:empirical-faithfulness}). Theorem~\ref{thm:commutation} then follows from the classical commutation criterion in finite dimension and from the Wasserstein stability of the nonlocal flows. We refer to \cite{PaulTrelat2026} for a survey on the passage from particle systems to mean-field and graph limits, and to \cite{PaulRossiTrelat2025} for systems with multiple-wise interactions.

\paragraph{Complete Lie algebras.}
Under Assumption~\ref{ass:regularity}, the bracket $\br{f}{g}$ is only continuous in general: it need not generate a flow, and iterated brackets are not defined. We therefore introduce in Section~\ref{sec:symbols} classes of nonlocal vector fields which are stable under brackets and whose elements generate global flows. The class $\SL$ consists of the polynomial nonlocal fields
$$
f(x,\mu)=\sum_{m=0}^M\int_{(\R^d)^m}K_m(x;y_1,\ldots,y_m)\,d\mu^{\otimes m}(y_1,\ldots,y_m),
$$
whose kernels satisfy global symbol estimates inspired by pseudodifferential calculus (see \cite{Hormander2007III}), and the class $\mathcal{B}_L$ consists of fields of the same form with bounded smooth kernels. By Corollary~\ref{cor:particle-morphisms}, the map $f\mapsto(f^N)_{N\geq1}$ is an injective Lie algebra morphism from $\SL$ to the product of the Lie algebras of vector fields on $(\R^d)^N$.

\paragraph{Quantitative formulas.}
For $f,g\in\SL$, we prove that
\begin{align*}
&\Phi_{-t}^g\circ\Phi_{-t}^f\circ\Phi_t^g\circ\Phi_t^f(\mu)=\Phi_{t^2}^{\br{f}{g}}(\mu)+\mathrm{O}\big((1+M_2(\mu))\vert t\vert^3\big),\\
&\Phi_t^g\circ\Phi_t^f(\mu)=\Phi_t^{f+g+\frac{t}{2}\br{f}{g}}(\mu)+\mathrm{O}\big((1+M_2(\mu))\vert t\vert^3\big),
\end{align*}
in $W_2$, where $M_2(\mu)=(\int\vert x\vert^2\,d\mu(x))^{1/2}$, with constants independent of $\mu$ and of the number of particles (Theorem~\ref{thm:quantitative}). The first expansion shows that the cycle of the four flows moves the population along the flow of the bracket during the time $t^2$. The second one is a second-order Baker-Campbell-Hausdorff (BCH) formula: applying $f$ and then $g$ is equivalent, up to the third order, to following the field $f+g$ corrected by $\frac{t}{2}\br{f}{g}$. The sign of this correction corresponds to applying $f$ before $g$; reversing the order replaces $\br{f}{g}$ by $-\br{f}{g}$. We prove these estimates in particle coordinates, with a norm weighted by the initial positions. The symbol estimates lead to constants independent of both the configuration and the number of particles, and passage to the limit requires only a finite second moment. In other words, the limit of fast switching and the mean-field limit commute at the second order. For bounded smooth kernels, the same argument gives a remainder $C\vert t\vert^3$, independent of the law (Proposition~\ref{prop:bounded-kernels}). Remark~\ref{rem:splitting} interprets the second expansion in terms of operator splitting.

\paragraph{Observables.}
We also study the action of these fields on scalar observables. The Lie-Lions derivatives $\mathscr{L}_fF(\mu)=\int\partial_\mu F(\mu)\cdot f(\cdot,\mu)\,d\mu$ satisfy $[\mathscr{L}_f,\mathscr{L}_g]=\mathscr{L}_{\br{f}{g}}$ and a Leibniz rule on an algebra $\mathcal{A}_L$ of polynomial observables, as vector fields and functions on a manifold. At a fixed law, the action of $f$ is represented by the orthogonal projection of $f(\cdot,\mu)$ onto the Wasserstein tangent space (Section~\ref{sec:observables}).

\paragraph{Examples.}
In Section~\ref{sec:examples}, we first show that fields which are equivariant under a rigid motion of the population commute with it, although their frozen spatial bracket does not vanish in general (Proposition~\ref{prop:symmetries}). We then consider models in which aggregation changes a collective rotation rate, coherence changes a common frequency of oscillators, a hierarchy transmits the turn of a leader with a delay increasing with its depth, two alternating pairwise interactions create an effective three-body interaction, and inward transport changes an occupancy-dependent circulation rate. The rotation model has explicit flows, the occupancy model belongs to $\SL$, and the other models have smooth bounded kernels.

\paragraph{Related work.}
Differential calculus on Wasserstein space can be developed through continuity equations and gradient tangent spaces (see \cite{AGS,DingFang2021,GangboKimPacini2011,Lott2008,Otto2001}), or through Lions differentiation of law-invariant lifts (see \cite{AlfonsiJourdain2020,Bertucci2026,Card12,CarmonaDelarue2018I,GangboTudorascu2019}).

The nonlocal correction in the bracket above already appears in several recent works. In \cite{DelloSchiavoPinamonti2026}, cylindrical derivations on probability measures over a closed manifold are endowed with such a bracket and a Lie-Rinehart structure (Corollary~2.13, formulated there as a Lie algebroid). In \cite{GoncharovaPogodaevStaritsyn2026}, nonlocal brackets defined through smooth displacement lifts are used to prove local projected controllability: \cite[Proposition~2 and Remark~5]{GoncharovaPogodaevStaritsyn2026} give the same decomposition into the spatial bracket and the measure correction.
Chapter~3 of \cite{Vicari2026} develops this bracket using Hadamard differentiability, with applications to variations of flows and pullbacks. In particular, \cite[Proposition~3.32]{Vicari2026} establishes a first-order pullback expansion with a quadratic remainder, uniformly on bounded subsets of $L^\infty$.

For locally Lipschitz vector fields which do not depend on the law, commutation of the flows is characterized in \cite{RampazzoSussmann2007} by means of set-valued Lie brackets, and a measure-theoretic version of these brackets is constructed in \cite{CavagnariMarigonda2018}. Here we establish the global commutation criterion under first-order assumptions on all of $\Ptwo$. We construct a complete Lie algebra for which the commutator and BCH estimates depend only on the second moment, with constants uniform in the particle number. The corresponding finite-dimensional expansions are classical (see \cite{AgrachevSachkov2004,BeauchardLeBorgneMarbach2023,MauhartMichor1992}). To our knowledge, Theorems~\ref{thm:commutation} and~\ref{thm:quantitative} do not appear as such in the literature. We learned of \cite{GoncharovaPogodaevStaritsyn2026,Vicari2026} while completing this work.

\medskip
The paper is organized as follows. Section~\ref{sec:bracket} is devoted to the first-order calculus, to particle realizations and to the proof of Theorem~\ref{thm:commutation}. In Section~\ref{sec:symbols}, we construct $\SL$ and $\mathcal{B}_L$, and we prove Theorem~\ref{thm:quantitative}. Section~\ref{sec:observables} concerns Lie-Lions derivatives of observables. The examples are developed in Section~\ref{sec:examples}. Section~\ref{sec:conclusion} gathers some perspectives. The appendix contains the differentiation and uniform particle estimates used in the proofs.

\section{Nonlocal flows and commutation}
\label{sec:bracket}

\subsection{Regularity and global flows}
\label{subsec:flows}

We set $\mathcal{H}=L^2((0,1);\R^d)$, where $(0,1)$ carries the Lebesgue measure $\mathbb{P}$. We write $\Vert \cdot\Vert_2$ for its norm, $\mathbb{E}$ for expectation, and $\mu_X=X_\#\mathbb{P}$ for the law of $X\in\mathcal{H}$. We identify the dual of $\mathcal{H}$ with $\mathcal{H}$ through the scalar product $\mathbb{E}[X\cdot Y]$. We recall that the Wasserstein distance of order $2$ is defined by (see \cite[Chapter~7]{AGS})
$$
W_2(\mu,\nu) =\Big(\inf_{\pi\in\Pi(\mu,\nu)} \int_{\R^d\times\R^d}\vert x-y\vert^2\,d\pi(x,y)\Big)^{1/2},
$$
where $\Pi(\mu,\nu)$ is the set of Borel probability measures on $\R^d\times\R^d$ with marginals $\mu$ and $\nu$, and we set $M_2(\mu)=(\int\vert x\vert^2\,d\mu(x))^{1/2}$.

\begin{definition}\label{def:Lions-derivative}
Let $k\geq1$ and let $F:\R^d\times\Ptwo\to\R^k$. We set $\widetilde F(x,X)=F(x,\mu_X)$ for $x\in\R^d$ and $X\in\mathcal{H}$. We say that $F$ is \emph{Lions-differentiable} at $(x,\mu)$ if $\widetilde F(x,\cdot)$ is Fr\'echet differentiable at some $X\in\mathcal{H}$ with law $\mu$. Its derivative is then represented by a unique $\partial_\mu F(x,\mu)\in L^2(\mu;\R^{k\times d})$ according to
\begin{equation}\label{lions-representation}
\partial_X\widetilde F(x,X)Y=\mathbb{E}\big[\partial_\mu F(x,\mu)(X)Y\big]
\end{equation}
for every $Y\in\mathcal{H}$ and for every $X\in\mathcal{H}$ with law $\mu$.
\end{definition}

\begin{remark}\label{rem:representation}
Let $F:\Ptwo\to\R$ be Lions-differentiable at $\mu$, and let $X\in\mathcal{H}$ with law $\mu$ be such that $\widetilde F$ is Fr\'echet differentiable at $X$. The Riesz theorem gives some $Z\in\mathcal{H}$ such that $D\widetilde F(X)Y=\mathbb{E}[Z\cdot Y]$ for every $Y\in\mathcal{H}$. Definition~\ref{def:Lions-derivative} asserts two further facts, which are not obvious: $Z$ is a function of $X$, namely $Z=\partial_\mu F(\mu)(X)$, and this function depends only on $\mu$, so that Fr\'echet differentiability at one random variable with law $\mu$ implies Fr\'echet differentiability at all of them, with the same representative. Both facts rely on law invariance, that is, $\widetilde F(X)=\widetilde F(X')$ whenever $X$ and $X'$ have the same law. However, two random variables with the same law are in general not obtained from each other by composition with a measure-preserving bijection of $(0,1)$, as shown by $X(\omega)=\omega$ and $X'(\omega)=2\omega$ modulo $1$, so that the chain rule does not apply directly. When $\widetilde F$ is continuously differentiable, which is the setting of Assumption~\ref{ass:regularity} below, we refer to \cite[Theorems~6.2 and~6.5]{Card12} and \cite[Propositions~5.24 and~5.25]{CarmonaDelarue2018I}. For a Fr\'echet derivative at a single point, we refer to \cite{AlfonsiJourdain2020} and \cite[Corollary~3.22]{GangboTudorascu2019}.
\end{remark}

A \emph{nonlocal vector field} is a map $f:\R^d\times\Ptwo\to\R^d$, and Definition~\ref{def:Lions-derivative} with $k=d$ gives $\partial_\mu f(x,\mu)(y)\in\R^{d\times d}$. A map defined on $\Ptwo$ is identified with a map independent of $x$. For $k=1$, we identify $\partial_\mu F(\mu)(y)\in\R^{1\times d}$ with a column vector of $\R^d$, so that \eqref{lions-representation} reads $D\widetilde F(X)Y=\mathbb{E}[\partial_\mu F(\mu)(X)\cdot Y]$. The \emph{Lagrangian lift} of $f$ is the map $\widehat f:\operatorname{Dom}(\widehat f)\to\mathcal{H}$ defined by $\widehat f(X)(\omega)=f(X(\omega),\mu_X)$ for almost every $\omega\in(0,1)$, where $\operatorname{Dom}(\widehat f)=\{X\in\mathcal{H}\ \mid\ f(X(\cdot),\mu_X)\in\mathcal{H}\}$. The map $\widetilde f$ is used to state the regularity of $f$, and $\widehat f$ gives the Lagrangian form $\dot X=\widehat f(X)$ of the continuity equation.

\begin{assumption}
\label{ass:regularity}
The map $\widetilde f:\R^d\times \mathcal{H}\to\R^d$ is continuously Fr\'echet differentiable. There is a constant $C>0$ such that
\begin{equation}\label{field-growth}
\vert f(x,\mu)\vert \leq C\big(1+\vert x\vert +M_2(\mu)\big) \qquad \forall x\in\R^d \qquad\forall \mu\in\Ptwo.
\end{equation}
For every $R>0$ there is $C_R>0$ such that, for $\Vert X\Vert_2\leq R$,
\begin{equation}\label{field-derivative-bounds}
\Vert \partial_x\widetilde f(x,X)\Vert \leq C_R, \qquad \Vert \partial_X\widetilde f(x,X)\Vert \leq C_R(1+\vert x\vert),
\end{equation}
where the second norm is the operator norm of linear maps from $\mathcal{H}$ to $\R^d$.
\end{assumption}

The growth bound \eqref{field-growth} ensures that $\operatorname{Dom}(\widehat f)=\mathcal{H}$. The second estimate in \eqref{field-derivative-bounds} is equivalent, up to a constant depending only on $d$, to $\Vert \partial_\mu f(x,\mu)\Vert_{L^2(\mu;\R^{d\times d})} \leq C_R(1+\vert x\vert)$ for $M_2(\mu)\leq R$, with the Hilbert-Schmidt norm on matrices. In particular, this assumption allows velocities depending on the barycenter, and on second moments through bounded smooth functions (see Section~\ref{subsec:skater}).

We recall the standard existence and stability result for nonlocal continuity equations (see \cite{CrippaLecureuxMercier2013,Dobrushin1979,PicRosNonlocalVelocity}). We give a short proof for the growth conditions of Assumption~\ref{ass:regularity}.

\begin{proposition}
\label{prop:global-flow}
Let $f$ satisfy Assumption~\ref{ass:regularity}. For every $\mu_0\in\Ptwo$, the continuity equation
\begin{equation}\label{nonlocal-continuity}
\partial_t\mu_t+\operatorname{div}\big(f(\cdot,\mu_t)\mu_t\big)=0, \qquad \mu_{t=0}=\mu_0,
\end{equation}
has a unique solution in $C(\R;\Ptwo)$, in the weak sense against $C_c^\infty(\R^d)$. It defines a flow $(\Phi_t^f)_{t\in\R}$ on $\Ptwo$. For every $R,T>0$ there are $R_T,L_{R,T}>0$ such that
\begin{equation}\label{W2-stability}
M_2(\Phi_t^f\mu)\leq R_T,\qquad
W_2(\Phi_t^f\mu,\Phi_t^f\nu) \leq e^{L_{R,T}\vert t\vert }W_2(\mu,\nu) \quad\text{for }\vert t\vert \leq T,
\end{equation}
whenever $M_2(\mu),M_2(\nu)\leq R$. Moreover, for every $s,t\in\R$,
\begin{equation}\label{nonlocal-flow-property}
\Phi_0^f=\operatorname{Id},\qquad \Phi_{t+s}^f=\Phi_t^f\circ\Phi_s^f,\qquad (\Phi_t^f)^{-1}=\Phi_{-t}^f.
\end{equation}
\end{proposition}

\begin{proof}
By \eqref{field-derivative-bounds}, for $\Vert X\Vert_2,\Vert Y\Vert_2\leq R$,
\begin{equation}\label{diagonal-local-Lipschitz}
\Vert\widehat f(X)-\widehat f(Y)\Vert_2\leq C_R(2+R)\Vert X-Y\Vert_2.
\end{equation}
Moreover, \eqref{field-growth} implies $\Vert\widehat f(X)\Vert_2\leq C(1+2\Vert X\Vert_2)$. The Cauchy-Lipschitz theorem in $\mathcal{H}$ and Gronwall's inequality therefore imply the existence of a unique global solution of $\dot X_t=\widehat f(X_t)$, with $\Vert X_t\Vert_2\leq(1+\Vert X_0\Vert_2)e^{2C\vert t\vert}-1$.

Let $\mu_t=\mu_{X_t}$, and let $T_t$ be the flow at time $t$ of the time-dependent field $f(\cdot,\mu_t)$, which is Lipschitz with respect to $x$, uniformly for $t$ in compact intervals, by \eqref{field-derivative-bounds}. Uniqueness of the spatial characteristics of $f(\cdot,\mu_t)$ implies $X_t=T_t(X_0)$ almost surely. For any $Y_0$ with the same law as $X_0$, the curve $T_t(Y_0)$ solves the same lifted equation and has law $\mu_t$. Hence $\mu_t$ depends only on $\mu_0$, and the chain rule gives \eqref{nonlocal-continuity}.

Conversely, every weak solution on a compact time interval has bounded second moments and square-integrable velocity. The superposition principle \cite[Theorem~8.2.1]{AGS}, with the resulting path measure realized on $(0,1)$, represents it as the law of an absolutely continuous solution of the same equation in $\mathcal{H}$. Uniqueness follows, including for negative times by replacing $f$ with $-f$, and implies \eqref{nonlocal-flow-property}.

Finally, we realize an optimal coupling of the initial laws in $\mathcal{H}\times \mathcal{H}$. Both solutions remain in the ball of radius $R_T=(1+R)e^{2CT}-1$. Applying Gronwall's inequality to \eqref{diagonal-local-Lipschitz} gives \eqref{W2-stability} with $L_{R,T}=C_{R_T}(2+R_T)$.
\end{proof}

For $R>0$, we set $B_R=\{\mu\in\Ptwo\ \mid\ M_2(\mu)\leq R\}$. We denote by $\operatorname{Aut}(\Ptwo)$ the group of homeomorphisms $\Psi:\Ptwo\to\Ptwo$ such that, for every $R>0$, there is $L_R>0$ with
\begin{equation}\label{automorphism-bounds}
\max\big( W_2(\Psi\mu,\Psi\nu), W_2(\Psi^{-1}\mu,\Psi^{-1}\nu) \big) \leq L_R W_2(\mu,\nu) \qquad\forall \mu,\nu\in B_R.
\end{equation}

These maps send moment balls into moment-bounded sets, since $M_2(\Psi\mu)\leq W_2(\Psi\mu,\Psi\delta_0)+M_2(\Psi\delta_0)$. Therefore the maps satisfying \eqref{automorphism-bounds} are stable under composition, and their inverses satisfy the same bounds by definition. By Proposition~\ref{prop:global-flow}, we have $\Phi_t^f\in\operatorname{Aut}(\Ptwo)$.

\begin{definition}\label{def:exponential}
For a field $f$ generating a unique global flow $(\Phi_t^f)_{t\in\R}$, we define its \emph{exponential} by $\exp(f)=\Phi_1^f$. We have $\exp(tf)=\Phi_t^f$ and $\exp((s+t)f)=\exp(sf)\circ\exp(tf)$.
\end{definition}

In particular, $\exp(0)=\Id$ and $\exp(f)^{-1}=\exp(-f)$.

\subsection{Lie-Lions bracket and particle realizations}\label{subsec:particles}

\begin{definition}\label{def:Lie-Lions-derivative}
Let $f$ be a nonlocal vector field such that $f(\cdot,\mu)\in L^2(\mu;\R^d)$ for every $\mu\in\Ptwo$, and let $F:\R^d\times\Ptwo\to\R^k$ be Lions-differentiable. The \emph{Lie-Lions derivative} of $F$ along $f$ is
\begin{equation}\label{Lie-Lions-derivative}
\mathscr{L}_fF(x,\mu)=\int_{\R^d}\partial_\mu F(x,\mu)(y)f(y,\mu)\,d\mu(y).
\end{equation}
\end{definition}

The integral converges by Cauchy-Schwarz. If $f$ satisfies Assumption~\ref{ass:regularity} and $\widetilde F(x,\cdot)$ is $C^1$, the chain rule in $\mathcal{H}$ applied to $\widetilde F(x,X_t)$, with $\dot X_t=\widehat f(X_t)$, gives $\frac{d}{dt}F(x,\Phi_t^f\mu)=\mathscr{L}_fF(x,\Phi_t^f\mu)$. Hence $\mathscr{L}_fF$ is the derivative of $F$ at fixed $x$ when the law is transported by $f$.

Let $f,g$ satisfy Assumption~\ref{ass:regularity}. The \emph{total derivative} of $g$ along $f$ is
\begin{equation}\label{total-derivative}
\mathcal{D}_fg=\partial_xg\,f+\mathscr{L}_fg.
\end{equation}
The first term differentiates $g$ along $f$ in the variable $x$, at fixed law, and the second one along the flow of $f$ in the variable $\mu$, at fixed $x$. Along $\mu_t=\Phi_t^f\mu$ and $\dot x_t=f(x_t,\mu_t)$, we have $\frac{d}{dt}g(x_t,\mu_t)=\mathcal{D}_fg(x_t,\mu_t)$.

\begin{definition}\label{def:Lie-Lions-bracket}
The \emph{Lie-Lions bracket} of $f$ and $g$ is
\begin{equation}\label{bracket-def}
\br{f}{g}=\mathcal{D}_fg-\mathcal{D}_gf=[f,g]_x+[f,g]_\mu,
\end{equation}
where $[f,g]_x=\partial_xg\,f-\partial_xf\,g$ and $[f,g]_\mu=\mathscr{L}_fg-\mathscr{L}_gf$.
\end{definition}

In terms of lifts, the total derivative is the Fr\'echet derivative of $\widetilde g$ in the direction $(\widetilde f,\widehat f)$:
\begin{equation}\label{lifted-total-derivative}
\mathcal{D}_fg(x,\mu_X)=D\widetilde g(x,X)\big(\widetilde f(x,X),\widehat f(X)\big).
\end{equation}
Since a $W_2$-convergent sequence of laws can be realized by an $\mathcal{H}$-convergent sequence, \eqref{lifted-total-derivative} implies continuity of $\mathcal{D}_fg$ and of $\br{f}{g}$ in $(x,\mu)$. Moreover, by Assumption~\ref{ass:regularity}, we have $\vert\br{f}{g}(x,\mu)\vert\leq K_R(1+\vert x\vert)$ whenever $M_2(\mu)\leq R$. In particular, $\br{f}{g}(\cdot,\mu)\in L^2(\mu)$.

\begin{definition}
\label{def:particle-realization}
For a nonlocal vector field $f$ and $N\geq1$, its \emph{particle realization} is the field $f^N=(f_i^N)_{1\leq i\leq N}$ defined, for $\boldsymbol{x}=(x_1,\ldots,x_N)\in(\R^d)^N$, by
$$
\mu_{\boldsymbol{x}}^N=\frac{1}{N}\sum_{i=1}^N\delta_{x_i}, \qquad f^N_i(\boldsymbol{x})=f(x_i,\mu_{\boldsymbol{x}}^N).
$$
\end{definition}

\begin{proposition}
\label{prop:particle}
Let $f,g$ satisfy Assumption~\ref{ass:regularity}. The particle fields $f^N,g^N$ are $C^1$ and complete. For $u\in\{f,g\}$, their flows solve
\begin{equation}\label{particle-ode}
\frac{d}{dt}\Phi_t^{u^N}(\boldsymbol x) =u^N\big(\Phi_t^{u^N}(\boldsymbol x)\big),\qquad \Phi_0^{u^N}(\boldsymbol x)=\boldsymbol x.
\end{equation}
Moreover,
\begin{equation}\label{empirical-flow}
\Phi_t^f\mu_{\boldsymbol{x}}^N =\mu_{\Phi_t^{f^N}(\boldsymbol{x})}^N.
\end{equation}
For all $\boldsymbol x,\boldsymbol h\in(\R^d)^N$, including configurations with coincident particle positions ($x_i=x_j$ for some $i\ne j$),
\begin{equation}\label{particle-derivative}
(Df^N(\boldsymbol{x})\boldsymbol{h})_i =\partial_xf(x_i,\mu_{\boldsymbol{x}}^N)h_i +\frac{1}{N}\sum_{j=1}^N \partial_\mu f(x_i,\mu_{\boldsymbol{x}}^N)(x_j)h_j.
\end{equation}
In particular, for every nonlocal vector field $u$ satisfying Assumption~\ref{ass:regularity}, we have
\begin{equation}\label{particle-total-derivative}
Df^N\,u^N=(\mathcal{D}_uf)^N,
\end{equation}
and
\begin{equation}\label{particle-bracket-compatibility}
[f^N,g^N]=\br{f}{g}^N.
\end{equation}
\end{proposition}

\begin{proof}
We partition $(0,1)$ into fixed intervals $I_1,\ldots,I_N$ of length $1/N$ and set $\mathcal{J}_N\boldsymbol q=\sum_{j=1}^Nq_j\mathbf{1}_{I_j}$. This defines a continuous linear map from $(\R^d)^N$ to $\mathcal{H}$, with $\mu_{\mathcal{J}_N\boldsymbol q}=\mu_{\boldsymbol q}^N$. Since $f_i^N(\boldsymbol q)=\widetilde f(q_i,\mathcal{J}_N\boldsymbol q)$, the particle field is $C^1$ and
$$
(Df^N(\boldsymbol q)\boldsymbol h)_i=\partial_x\widetilde f(q_i,\mathcal{J}_N\boldsymbol q)h_i+\partial_X\widetilde f(q_i,\mathcal{J}_N\boldsymbol q)[\mathcal{J}_N\boldsymbol h].
$$
Formula \eqref{lions-representation} identifies the last term with $N^{-1}\sum_{j=1}^N\partial_\mu f(q_i,\mu_{\boldsymbol q}^N)(q_j)h_j$, proving \eqref{particle-derivative} (see also \cite{DelarueLackerRamanan2019}). The fixed partition allows distinct perturbations of particles at the same position, so the formula also holds at collisions. Taking $\boldsymbol h=u^N(\boldsymbol q)$ in \eqref{particle-derivative}, the last term is $\mathscr{L}_uf(q_i,\mu_{\boldsymbol q}^N)$ by \eqref{Lie-Lions-derivative}, and \eqref{total-derivative} gives \eqref{particle-total-derivative}.

With $\Vert\boldsymbol q\Vert_N=\Vert\mathcal{J}_N\boldsymbol q\Vert_2$, \eqref{field-growth} implies $\Vert f^N(\boldsymbol q)\Vert_N\leq C(1+2\Vert\boldsymbol q\Vert_N)$. Hence $f^N$ is complete. The identity $\widehat f(\mathcal{J}_N\boldsymbol q)=\mathcal{J}_Nf^N(\boldsymbol q)$ and uniqueness of the lifted equation imply \eqref{empirical-flow}. Finally, \eqref{particle-total-derivative} applied to $g$ with $u=f$, and to $f$ with $u=g$, gives $[f^N,g^N]=Dg^N\,f^N-Df^N\,g^N=(\mathcal{D}_fg-\mathcal{D}_gf)^N$, that is, \eqref{particle-bracket-compatibility} by \eqref{bracket-def}.
\end{proof}

Identity \eqref{particle-bracket-compatibility} transfers the finite-dimensional Lie calculus to nonlocal fields, for every $N$, and all the proofs below rely on it. We shall also use the following consequence of the density of empirical measures.

\begin{lemma}
\label{lem:empirical-faithfulness}
Let $v:\R^d\times\Ptwo\to\R^d$ be continuous. If $v^N=0$ for every $N\geq1$, then $v=0$.
\end{lemma}

\begin{proof}
Let $(x,\mu)\in\R^d\times\Ptwo$ be arbitrary. By the density of finitely supported probability measures with rational weights \cite[proof of Proposition~7.1.5]{AGS}, we choose, along a sequence of integers $N\to+\infty$, equal-weight empirical measures $\nu_{N-1}$ with $N-1$ atoms, allowing repetitions, such that $W_2(\nu_{N-1},\mu)\to0$. Then $\mu_N=\frac{1}{N}\delta_x+\frac{N-1}{N}\nu_{N-1} \longrightarrow\mu$ in $W_2$, since $W_2(\mu_N,\mu)^2\leq N^{-1}W_2(\delta_x,\mu)^2+(1-N^{-1})W_2(\nu_{N-1},\mu)^2$. The point $x$ is one of the particles of $\mu_N$, so $v(x,\mu_N)=0$. By continuity, we get $v(x,\mu)=0$.
\end{proof}

\subsection{Commutation criterion}\label{subsec:commutation}

\begin{theorem}\label{thm:commutation}
Let $f,g$ satisfy Assumption~\ref{ass:regularity}. Their flows $\Phi_t^f,\Phi_t^g$ are defined on $\Ptwo$ for every $t\in\R$. The following assertions are equivalent:
\begin{enumerate}[label=\textnormal{(\roman*)}]
\item $\br{f}{g}(x,\mu)=0$ for every $x\in\R^d$ and $\mu\in\Ptwo$;
\item $\Phi_t^f\Phi_s^g\mu=\Phi_s^g\Phi_t^f\mu$ for every $s,t\in\R$ and $\mu\in\Ptwo$.
\end{enumerate}
\end{theorem}

\begin{proof}
Let us first assume that $\br{f}{g}=0$. By Proposition~\ref{prop:particle}, we have $[f^N,g^N]=0$ for every $N$. The complete $C^1$ particle fields therefore have commuting flows. Let us recall the classical argument: the identity $[f^N,g^N]=0$ implies $D\Phi_t^{f^N}(\boldsymbol{x})g^N(\boldsymbol{x}) =g^N(\Phi_t^{f^N}(\boldsymbol{x}))$ by the variational equation along $\Phi_t^{f^N}$. Uniqueness in \eqref{particle-ode} with $u=g$ then gives $\Phi_t^{f^N}\Phi_s^{g^N} =\Phi_s^{g^N}\Phi_t^{f^N}$ for all $s,t$. Let $\mu\in\Ptwo$. We choose equal-weight empirical laws $\mu_N\to\mu$ in $W_2$. By \eqref{empirical-flow}, $\Phi_t^f\Phi_s^g\mu_N=\Phi_s^g\Phi_t^f\mu_N$. Their initial second moments are bounded. For fixed $s,t$, the moment and stability estimates of Proposition~\ref{prop:global-flow} therefore give a constant $C_{s,t}$, independent of $N$, such that
\begin{align*}
&W_2(\Phi_t^f\Phi_s^g\mu_N,\Phi_t^f\Phi_s^g\mu) +W_2(\Phi_s^g\Phi_t^f\mu_N,\Phi_s^g\Phi_t^f\mu) \leq C_{s,t}W_2(\mu_N,\mu),\\
&W_2(\Phi_t^f\Phi_s^g\mu,\Phi_s^g\Phi_t^f\mu) \leq C_{s,t}W_2(\mu_N,\mu)\rightarrow0.
\end{align*}

Conversely, let us assume that the law flows commute for all initial laws and all times. Let $\boldsymbol{x}$ be a configuration with pairwise distinct components. By \eqref{empirical-flow}, we have $\mu_{\Phi_t^{f^N}\Phi_s^{g^N}(\boldsymbol{x})}^N =\mu_{\Phi_s^{g^N}\Phi_t^{f^N}(\boldsymbol{x})}^N$. For sufficiently small $\vert s\vert ,\vert t\vert $, the $i$th component of either configuration remains in a ball about $x_i$, with these balls pairwise disjoint. Equality of the empirical laws then implies equality of the labelled configurations. We differentiate this local commutation identity first in $s$ at zero and then in $t$ at zero to obtain $[f^N,g^N](\boldsymbol{x})=0$. Configurations with distinct components are dense, so, by continuity, we have $[f^N,g^N]=0$ everywhere. Proposition~\ref{prop:particle} and Lemma~\ref{lem:empirical-faithfulness}, applied to the continuous field $\br{f}{g}$, imply $\br{f}{g}=0$.
\end{proof}

\begin{remark}\label{rem:local-commutation}
The proof shows that both assertions of Theorem~\ref{thm:commutation} can be weakened. In (ii), it suffices that, for every equal-weight empirical law $\mu$ with distinct atoms, the identity $\Phi_t^f\Phi_s^g\mu=\Phi_s^g\Phi_t^f\mu$ holds for $\vert s\vert$ and $\vert t\vert$ small enough, depending on $\mu$. In (i), it suffices that $\br{f}{g}(\cdot,\mu)=0$ $\mu$-almost everywhere, for every $\mu\in\Ptwo$. Indeed, since the atoms of an empirical law have positive mass, this implies $\br{f}{g}^N=0$ for every $N$, and Lemma~\ref{lem:empirical-faithfulness} applies. In particular, commutation for small times on empirical laws implies commutation for all times on $\Ptwo$.
\end{remark}

\section{Complete Lie algebras and flow formulas}\label{sec:symbols}

Under Assumption~\ref{ass:regularity}, the bracket $\br{f}{g}$ is only continuous in general. It need not generate a flow, and iterated brackets are not defined. However, the estimates of Theorem~\ref{thm:quantitative} below involve the flows of $\br{f}{g}$ and of $Z_t=t(f+g)+\frac{t^2}{2}\br{f}{g}$, and a Lie algebra structure requires iterated brackets. We therefore construct \emph{complete} Lie algebras of nonlocal vector fields, that is, vector spaces stable under $\br{\cdot}{\cdot}$ whose elements generate unique global flows in $\operatorname{Aut}(\Ptwo)$.

The class $\SL$ of Section~\ref{subsec:polynomial-symbols} consists of fields that are polynomial in the law, $f(x,\mu)=\int K(x;y_1,\ldots,y_m)\,d\mu^{\otimes m}(y_1,\ldots,y_m)$, with kernels satisfying global symbol estimates. For $m=1$ and $m=2$, these fields are the pairwise and three-body interactions of mean-field models. The bracket preserves this structure: the Lions derivative lowers the degree in $\mu$ by one, and the composition in \eqref{total-derivative-kernel} adds the degrees (Proposition~\ref{prop:symbol}). Corollary~\ref{cor:particle-morphisms} realizes $\SL$, for each $N$, as a Lie algebra of smooth vector fields on $(\R^d)^N$, and Theorem~\ref{thm:quantitative} is proved in these finite-dimensional realizations (Section~\ref{subsec:quantitative}). The class $\mathcal{B}_L$ of Section~\ref{subsec:bounded-kernels} replaces the symbol estimates by bounded derivatives, and gives remainders independent of the law.

\subsection{Polynomial nonlocal symbols}\label{subsec:polynomial-symbols}

We set $\langle x\rangle=(1+\vert x\vert^2)^{1/2}$. For $m\geq0$, let $\mathcal{K}_m$ consist of the smooth maps $K:\R^d\times(\R^d)^m\longrightarrow\R^d$ which are symmetric in their last $m$ variables and satisfy, for all multi-indices $\alpha,\beta_1,\ldots,\beta_m$,
\begin{equation}\label{kernel-symbol}
\big\vert \partial_x^\alpha \partial_{y_1}^{\beta_1}\cdots \partial_{y_m}^{\beta_m}K(x;y_1,\ldots,y_m)\big\vert \leq C_{\alpha,\boldsymbol\beta} \langle x\rangle^{1-\vert \alpha\vert } \prod_{j=1}^m\langle y_j\rangle^{-\vert \beta_j\vert }.
\end{equation}
These estimates are inspired by symbol classes in pseudodifferential calculus. For $m=0$, they mean that each component of $K$ belongs to the class $S^1(\R^d)$ of symbols of order one (see \cite[Section~18.1]{Hormander2007III}). For $m\geq1$, they have the structure of the SG classes of Parenti \cite{Parenti1972} and Cordes \cite{Cordes1995}, where G stands for global: a symbol $a(x,\xi)$ of order $(m_1,m_2)$ satisfies $\vert\partial_x^\alpha\partial_\xi^\beta a(x,\xi)\vert\leq C_{\alpha,\beta}\langle x\rangle^{m_1-\vert\alpha\vert}\langle\xi\rangle^{m_2-\vert\beta\vert}$ on $\R^d\times\R^d$ (see \cite[Chapter~3]{NicolaRodino2010}, where they are called G-classes). In \eqref{kernel-symbol}, the variables $y_1,\ldots,y_m$ play the role of $\xi$, with order one in $x$ and order zero in each $y_j$. We use these global estimates for three reasons. First, they give fields with linear growth in $x$, uniformly with respect to the law, and thus global flows on $\Ptwo$. Second, the weights cancel in the nonlocal term of the bracket (see \eqref{kernel-weight-cancellation}), which makes the class stable under brackets. Third, they give bounds for the particle fields that are uniform in the configuration and in the number of particles (Lemma~\ref{lem:weighted-particles}). Estimates on compact sets only would not give global flows, and the constants would depend on the supports of the laws. Note that we only use these estimates, and not the pseudodifferential calculus. Symmetrization in $y_1,\ldots,y_m$ preserves both \eqref{kernel-symbol} and the integral against $\mu^{\otimes m}$.

\begin{definition}\label{def:polynomial-symbols}
For $m\geq0$, let $\SLm{m}$ be the set of fields of the form
\begin{equation}\label{polynomial-field}
f(x,\mu)=\int_{(\R^d)^m}K(x;y_1,\ldots,y_m)\, d\mu^{\otimes m}(y_1,\ldots,y_m), \qquad K\in\mathcal{K}_m.
\end{equation}
For $m=0$, this means $f(x,\mu)=K(x)$. The polynomial nonlocal symbol class is $\displaystyle \SL=\bigcup_{m\geq0}\SLm{m}$.
\end{definition}

These spaces form a filtration. Indeed, a kernel in $\mathcal{K}_m$ can be viewed as a kernel with an unused $(m+1)$st variable and then symmetrized. Its integral is unchanged because $\mu$ has total mass one. Therefore $\SLm{m}\subset \SLm{m+1}$, and $\SL$ is a vector space. Equivalently, its elements are finite sums of expressions \eqref{polynomial-field}. Polynomial functionals of measures and their degree filtrations are standard objects \cite{CuchieroLarssonSvalutoFerro2019}. On probability measures there is no intrinsic decomposition into homogeneous components of different degrees.

\begin{proposition}\label{prop:symbol}
The space $\SL$ is a Lie algebra under $\br{\cdot}{\cdot}$, and $\br{\SLm{m}}{\SLm{n}} \subset \SLm{m+n}$. Every $f\in\SL$ satisfies Assumption~\ref{ass:regularity} and generates a global flow on $\Ptwo$.
\end{proposition}

\begin{proof}
For $f$ given by \eqref{polynomial-field} with a kernel $K_f\in\mathcal{K}_m$, differentiation under the integral gives
\begin{align}
\notag \partial_xf(x,\mu) &=\int \partial_xK_f(x;\boldsymbol y)\,d\mu^{\otimes m}(\boldsymbol y),\\
\label{symbol-lions-derivative} \partial_\mu f(x,\mu)(z) &=m\int \partial_{y_1}K_f(x;z,y_2,\ldots,y_m)\, d\mu^{\otimes(m-1)}(y_2,\ldots,y_m)
\end{align}
when $m\geq1$, and $\partial_\mu f=0$ when $m=0$. In particular,
\begin{equation}\label{symbol-first-bounds}
\vert f(x,\mu)\vert \leq C\langle x\rangle,\qquad \vert \partial_xf(x,\mu)\vert \leq C,\qquad \vert \partial_\mu f(x,\mu)(z)\vert \leq C\langle x\rangle\langle z\rangle^{-1}.
\end{equation}
All constants are independent of $\mu$.

Let us justify the Fr\'echet regularity of $\widetilde f$. On the product probability space $(0,1)^m$, we have $\widetilde f(x,X)=\mathbb{E}\big[K_f(x;X(\omega_1),\ldots,X(\omega_m))\big]$. On each bounded set of positions $x$, the first and second derivatives of $K_f$ are bounded uniformly in the measure variables, by \eqref{kernel-symbol}. Taylor's formula and Cauchy-Schwarz give a remainder bounded by $C(\vert\xi\vert+\Vert Y\Vert_2)^2$ when $(x,X)$ is replaced by $(x+\xi,X+Y)$. The same derivative bounds imply continuity of the first derivative in operator norm. Hence $\widetilde f$ is $C^1$, with the derivatives above. Together with \eqref{symbol-first-bounds}, this proves Assumption~\ref{ass:regularity}, and Proposition~\ref{prop:global-flow} implies completeness.

Let $f$ and $g$ be represented by $K_f\in\mathcal{K}_m$ and $K_g\in\mathcal{K}_n$. We write $\boldsymbol y=(y_1,\ldots,y_m)$ and $\boldsymbol z=(z_1,\ldots,z_n)$. If $\operatorname{Sym}$ denotes averaging over all $m+n$ measure variables, a kernel for the total derivative \eqref{total-derivative} is
\begin{equation}\label{total-derivative-kernel}
T(K_f,K_g)(x;\boldsymbol y,\boldsymbol z) = \operatorname{Sym}\big[ \partial_xK_g(x;\boldsymbol z)K_f(x;\boldsymbol y) +n\partial_{z_1}K_g(x;\boldsymbol z)K_f(z_1;\boldsymbol y) \big].
\end{equation}
The second term is omitted when $n=0$. Formula \eqref{total-derivative-kernel} follows from \eqref{symbol-lions-derivative} and Fubini's theorem. Hence a kernel for $\br{f}{g}$ is $T(K_f,K_g)-T(K_g,K_f)$, where the variables in the second term are relabelled in the same list.

The local product $\partial_xK_g(x;\boldsymbol z)K_f(x;\boldsymbol y)$ satisfies \eqref{kernel-symbol} by Leibniz' rule. In the nonlocal product, $z_1$ is both a differentiated measure variable of $K_g$ and the spatial variable of $K_f$. At order zero the relevant weights cancel:
\begin{equation}\label{kernel-weight-cancellation}
\vert \partial_{z_1}K_g(x;\boldsymbol z)K_f(z_1;\boldsymbol y)\vert \leq C\frac{\langle x\rangle}{\langle z_1\rangle} \langle z_1\rangle =C\langle x\rangle.
\end{equation}
After $r$ further derivatives in $z_1$, a Leibniz term with $a$ derivatives on the first factor and $r-a$ on the second has weight $\langle z_1\rangle^{-1-a} \langle z_1\rangle^{1-(r-a)} =\langle z_1\rangle^{-r}$. Derivatives in $x$ act on $K_g$ only. Derivatives in the remaining measure variables have exactly their required negative weights. These observations, applied simultaneously to all multi-indices, show that $T(K_f,K_g)\in\mathcal{K}_{m+n}$. Hence both the total derivative and the bracket preserve the class, with $\mathcal{D}_fg,\ \br{f}{g}\in \SLm{m+n}$.

Bilinearity and antisymmetry are immediate. To prove the Jacobi identity, we use Proposition~\ref{prop:particle}: the particle realization of $J(f,g,h)=\br{\br{f}{g}}{h} +\br{\br{g}{h}}{f} +\br{\br{h}{f}}{g}$ is the ordinary Jacobi expression for three smooth particle fields, and therefore vanishes for every $N$.

Lemma~\ref{lem:empirical-faithfulness}, applied to $J(f,g,h)\in\SL$, proves the Jacobi identity.
\end{proof}

For a smooth manifold $M$, $\mathfrak{X}(M)$ denotes the Lie algebra of smooth vector fields on $M$.

\begin{corollary}\label{cor:particle-morphisms}
For every $N\geq1$, the map $\mathcal{R}_N:\SL\rightarrow \mathfrak{X}((\R^d)^N)$ defined by $\mathcal{R}_Nf=f^N$ is a Lie algebra morphism. The family $(\mathcal{R}_N)_{N\geq1}$ is jointly faithful, so that the map $\SL\rightarrow\prod_{N\geq1}\mathfrak{X}((\R^d)^N)$ given by $f\mapsto(f^N)_{N\geq1}$ is an injective Lie algebra morphism.
\end{corollary}

\begin{proof}
The symbol kernels make $f^N$ smooth. Linearity is immediate, bracket preservation follows from Proposition~\ref{prop:particle}, and injectivity from Lemma~\ref{lem:empirical-faithfulness}.
\end{proof}

\subsection{Commutator and BCH estimates}\label{subsec:quantitative}

Products of maps act from right to left. For $f,g\in\SL$, we set
$$
\mathcal{C}_t^{f,g}=\exp(-tg)\circ\exp(-tf)\circ\exp(tg)\circ\exp(tf),\qquad Z_t=t(f+g)+\frac{t^2}{2}\br{f}{g}.
$$
By Proposition~\ref{prop:symbol}, the fields $\br{f}{g}$ and $Z_t$ belong to $\SL$, and thus their exponentials are defined.

\begin{theorem}\label{thm:quantitative}
Let $f,g\in\SL$. There exist $C,t_0>0$ such that, for every $\mu\in\Ptwo$ and $\vert t\vert\leq t_0$,
\begin{align*}
W_2\big(\mathcal{C}_t^{f,g}\mu,\exp(t^2\br{f}{g})\mu\big)&\leq C\big(1+M_2(\mu)\big)\vert t\vert^3,\\
W_2\big(\exp(tg)\exp(tf)\mu,\exp(Z_t)\mu\big)&\leq C\big(1+M_2(\mu)\big)\vert t\vert^3.
\end{align*}
The constants depend on finitely many seminorms of kernel representations of $f,g$ and their measure degrees, and are independent of $\mu$.
\end{theorem}

The second estimate is a second-order Baker-Campbell-Hausdorff (BCH) formula. The sign of its bracket term corresponds to applying $f$ before $g$; reversing the order replaces $\br{f}{g}$ by $-\br{f}{g}$. We prove the estimates by applying finite-dimensional flow expansions in particle coordinates. Direct differentiation of the Lagrangian lift on $\mathcal{H}$ would require additional regularity, as the following example shows.

\begin{remark}\label{rem:diagonal-smoothness}
For the smooth field $f(x)=\sin x$, the Lagrangian lift $X\mapsto\sin X$ is not Fr\'echet differentiable at zero in $L^2$. Its directional derivatives at zero are those of the identity, but $X_\varepsilon=\pi\mathbf{1}_{(0,\varepsilon)}$ satisfies $\sin X_\varepsilon=0$ and $\Vert\sin X_\varepsilon-X_\varepsilon\Vert_2/\Vert X_\varepsilon\Vert_2=1$, while $X_\varepsilon\to0$ in $L^2$. This is a standard obstruction for superposition operators (see \cite[Section~3.6]{AppellZabrejko1990}).
\end{remark}

The particle realizations are restrictions of the Lagrangian lift. Indeed, the identity $\widehat f(\mathcal{J}_N\boldsymbol q)=\mathcal{J}_Nf^N(\boldsymbol q)$ used in the proof of Proposition~\ref{prop:particle} shows that the finite-dimensional subspace $\mathcal{J}_N((\R^d)^N)$ of $\mathcal{H}$ is invariant under $\widehat f$, and that, through the linear isomorphism $\mathcal{J}_N$ onto this subspace, $f^N$ is the restriction of $\widehat f$. This restriction is $C^1$ under Assumption~\ref{ass:regularity}. By \eqref{particle-bracket-compatibility}, the restriction of the Lagrangian lift of $\br{f}{g}$ is the Lie bracket of the restrictions of $\widehat f$ and $\widehat g$. The union of these subspaces over $N\geq1$ is dense in $\mathcal{H}$, and our proofs only differentiate these restrictions. In \cite[Proposition~3.20]{Vicari2026}, the lift of the bracket is identified with the commutator of the lifted fields, computed with Hadamard derivatives at bounded random variables.

\begin{proof}[Proof of Theorem~\ref{thm:quantitative}]
Let $\boldsymbol x=(x_1,\ldots,x_N)$ be an initial configuration. We set
$$
a_i=\langle x_i\rangle,\qquad \Vert \boldsymbol h\Vert_{\boldsymbol x,*} =\max_{1\leq i\leq N}\frac{\vert h_i\vert }{a_i}.
$$
Lemma~\ref{lem:weighted-particles} bounds $f^N,g^N$ and their first three derivatives on a fixed ball about $\boldsymbol x$ in this norm. Applying Lemma~\ref{lem:finite-dimensional-formulas} and Proposition~\ref{prop:particle}, we obtain both particle estimates with remainder $C\vert t\vert^3$, for $\vert t\vert\leq t_0$, where $C,t_0$ are independent of $N$ and $\boldsymbol x$.

If $\boldsymbol y,\boldsymbol z$ form either pair of resulting particle configurations, coupling equally labelled particles gives
$$
W_2\Big(\frac{1}{N}\sum_i\delta_{y_i}, \frac{1}{N}\sum_i\delta_{z_i}\Big) \leq \Big(\frac{1}{N}\sum_i\vert y_i-z_i\vert^2\Big)^{1/2} \leq C\vert t\vert^3 \Big(\frac{1}{N}\sum_i\langle x_i\rangle^2\Big)^{1/2}.
$$
The last factor is $(1+M_2(\mu_{\boldsymbol x}^N)^2)^{1/2}\leq1+M_2(\mu_{\boldsymbol x}^N)$. The particle realization of the field $Z_t$ is the corresponding finite-dimensional BCH field by linearity and Proposition~\ref{prop:particle}. This proves both estimates for empirical measures.

Let us now take equal-weight empirical measures $\mu_k\to\mu$ in $W_2$. Their second moments converge. For each fixed $t$, Proposition~\ref{prop:global-flow} gives continuity of all the flow maps in the two estimates, including those of $\br{f}{g}$ and $Z_t$, since these fields belong to $\SL$. Passing to the limit proves both estimates for $\mu$.
\end{proof}

\begin{remark}\label{rem:splitting}
The second estimate of Theorem~\ref{thm:quantitative} concerns the Lie splitting, which approximates the flow of $f+g$ by the composition $\exp(tg)\circ\exp(tf)$. Since $\exp(Z_t)=\Phi_t^{f+g+\frac{t}{2}\br{f}{g}}$, this composition coincides with the flow of $f+g$ corrected by $\frac{t}{2}\br{f}{g}$, up to a remainder of order three in $W_2$, uniformly in the number of particles. If $f$ and $g$ satisfy Assumption~\ref{ass:regularity} and $\br{f}{g}=0$, then the splitting is exact, i.e., $\exp(tg)\circ\exp(tf)=\exp(t(f+g))$ for every $t\in\R$. Indeed, the commuting particle fields satisfy $\Phi_t^{g^N}\circ\Phi_t^{f^N}=\Phi_t^{f^N+g^N}$, and we conclude by \eqref{empirical-flow} and by density of equal-weight empirical laws, as in the proof of Theorem~\ref{thm:commutation}.
\end{remark}

\subsection{Smooth bounded kernels}\label{subsec:bounded-kernels}
The proof of the estimates of Theorem~\ref{thm:quantitative} in Section~\ref{subsec:quantitative} uses the symbol estimates only to bound the particle fields and their first three derivatives, uniformly in the number of particles and in the configuration (Lemma~\ref{lem:weighted-particles}). For kernels with bounded derivatives, such bounds hold for the unweighted norm $\max_i\vert h_i\vert$, and the factor $1+M_2(\mu)$ disappears from the remainder. Let $\mathcal{B}_L$ denote the set of fields of the form \eqref{polynomial-field}, with a symmetric kernel $K \in C_b^\infty(\R^d \times (\R^d)^m;\R^d)$, $m \geq 0$, where $C_b^\infty$ means that the kernel and all its derivatives are bounded. As in Definition~\ref{def:polynomial-symbols}, finite sums are included by adding unused integration variables.

The fields of $\mathcal{B}_L$ have the same polynomial form as those of $\SL$, and only the conditions on the kernels differ. In $\SL$, a kernel may grow linearly in $x$, and each derivative brings a decay factor $\langle x\rangle^{-1}$ or $\langle y_j\rangle^{-1}$. In $\mathcal{B}_L$, a kernel and all its derivatives are bounded, without decay. Neither class contains the other. The linear field $f(x,\mu)=Ax$, with $A\neq0$, belongs to $\SL$ but not to $\mathcal{B}_L$, since the fields of $\mathcal{B}_L$ are bounded. Conversely, the convolution fields $f(x,\mu)=\int\psi(x-y)\,d\mu(y)$ with $\psi\in C_b^\infty(\R^d;\R^d)$ belong to $\mathcal{B}_L$, and to $\SL$ only when $\psi$ is constant. Indeed, $\partial_\mu f(x,\mu)(z)=-D\psi(x-z)$. Taking $\mu$ with full support, if $f\in\SL$, then \eqref{kernel-symbol} and \eqref{symbol-lions-derivative} imply $\vert D\psi(x-z)\vert\leq C\langle x\rangle\langle z\rangle^{-1}$ and $\vert D^2\psi(x-z)\vert\leq C\langle z\rangle^{-1}$ for all $x,z\in\R^d$. Letting $z$ tend to infinity with $x-z$ fixed, we get $D^2\psi=0$, and then, with $x$ fixed, $D\psi=0$. The class $\mathcal{B}_L$ also contains the periodic heading kernels used in Section~\ref{subsec:flock-turning}. On the other hand, $\SL$ contains the fields $\sum_k a_k(x)F_k(\mu)$, where each $a_k$ satisfies \eqref{kernel-symbol} with $m=0$ and each $F_k$ belongs to the algebra $\mathcal{A}_L$ of Section~\ref{subsec:anchor}, such as the fields of Section~\ref{subsec:congestion}.

\begin{proposition}\label{prop:bounded-kernels}
The space $\mathcal{B}_L$ is a complete Lie algebra under $\br{\cdot}{\cdot}$, and every field in $\mathcal{B}_L$ satisfies Assumption~\ref{ass:regularity}. For $f,g \in \mathcal{B}_L$, there exist $C,t_0>0$ such that, for every $\mu \in \Ptwo$ and $\vert t\vert \leq t_0$,
\begin{align*}
W_2\big(\mathcal{C}_t^{f,g}\mu,\exp(t^2\br{f}{g})\mu\big) &\leq C\vert t\vert^3,\\
W_2\big(\exp(tg)\exp(tf)\mu,\exp(Z_t)\mu\big) &\leq C\vert t\vert^3,
\end{align*}
where $\mathcal{C}_t^{f,g}$ and $Z_t=t(f+g)+\frac{t^2}{2}\br{f}{g}$ are defined as in Section~\ref{subsec:quantitative}. The constants depend only on $d$, the measure degrees and the bounds for derivatives of order at most three of chosen kernel representations of $f,g$. They are independent of $\mu$ and of the number of particles.
\end{proposition}

\begin{proof}
The map $\widetilde f$ is $C^1$ by the Taylor argument used in the proof of Proposition~\ref{prop:symbol}. Here $f$, $\partial_xf$ and $\partial_X\widetilde f$ are bounded, and thus Assumption~\ref{ass:regularity} holds. Completeness follows from Proposition~\ref{prop:global-flow}. Formula \eqref{total-derivative-kernel} proves closure under brackets, since products and derivatives of $C_b^\infty$ kernels remain in $C_b^\infty$. The Jacobi identity follows from particle compatibility and Lemma~\ref{lem:empirical-faithfulness}.

For the quantitative estimates, we equip $(\R^d)^N$ with the norm $\Vert\boldsymbol h\Vert_\infty=\max_i\vert h_i\vert$. If $f$ has a kernel $K$ with $m$ measure variables, then
$$
f_i^N(\boldsymbol y)=\frac{1}{N^m}\sum_{j_1,\ldots,j_m=1}^N K(y_i;y_{j_1},\ldots,y_{j_m}).
$$
A derivative of order $j$ is a sum of terms obtained by assigning each of the $j$ differentiations to one of the $m+1$ slots. Bounded kernel derivatives therefore give
$$
\Vert D^j f^N(\boldsymbol y)[\boldsymbol h^1,\ldots,\boldsymbol h^j]\Vert_\infty \leq C_j\prod_{\ell=1}^j\Vert\boldsymbol h^\ell\Vert_\infty, \qquad 0 \leq j \leq 3,
$$
uniformly in $N$ and $\boldsymbol y$, with the empty product equal to one. The normalization by $N^m$ cancels the number of summands, including when indices coincide, since differentiation is performed on the kernel slots. The same estimates hold for $g^N$.

We apply Lemma~\ref{lem:finite-dimensional-formulas} on balls of radius one for this norm, with constants independent of their centers and of $N$. By Proposition~\ref{prop:particle}, the resulting particle fields for the commutator and the BCH approximation are $(\br{f}{g})^N$ and $Z_t^N$. Both errors are bounded by $C\vert t\vert^3$ in $\Vert\cdot\Vert_\infty$. For either pair of final configurations $\boldsymbol y,\boldsymbol z$, the coupling of equally labelled particles gives
$$
W_2(\mu_{\boldsymbol y}^N,\mu_{\boldsymbol z}^N) \leq \Big(\frac{1}{N}\sum_{i=1}^N\vert y_i-z_i\vert^2\Big)^{1/2} \leq \Vert\boldsymbol y-\boldsymbol z\Vert_\infty \leq C\vert t\vert^3.
$$
Finally, equal-weight empirical measures are dense in $\Ptwo$. Since $f,g,\br{f}{g},Z_t \in \mathcal{B}_L$, their flows are continuous in $W_2$ by Proposition~\ref{prop:global-flow}. Passing to the limit proves the two estimates for every $\mu \in \Ptwo$, with the same constants.
\end{proof}

\section{Differential operators on observables}\label{sec:observables}

\subsection{Lie-Lions derivatives}\label{subsec:Lie-derivatives}

We now view the Lie-Lions derivatives \eqref{Lie-Lions-derivative} as operators on scalar observables. For $F_\varphi(\mu)=\int \varphi\,d\mu$ with $\varphi\in C_c^\infty(\R^d)$, we have $\mathscr{L}_fF_\varphi(\mu)=\int\nabla \varphi(x)\cdot f(x,\mu)\,d\mu(x)$.

\begin{proposition}\label{prop:operator}
If $f,g$ satisfy Assumption~\ref{ass:regularity}, then
\begin{equation}\label{operator-bracket}
(\mathscr{L}_f\mathscr{L}_g-\mathscr{L}_g\mathscr{L}_f)F_\varphi =\mathscr{L}_{\br{f}{g}}F_\varphi \qquad \forall \varphi\in C_c^\infty(\R^d).
\end{equation}
\end{proposition}

\begin{proof}
By Lemma~\ref{lem:iterability} in Appendix~\ref{app:calculus}, applied to $g$, $\mathscr{L}_gF_\varphi$ is Lions-differentiable and we have
$$
\mathscr{L}_f\mathscr{L}_gF_\varphi(\mu) = \int f(x,\mu)^TD^2\varphi(x)g(x,\mu)\,d\mu(x) + \int\nabla \varphi(x)\cdot \mathcal{D}_fg(x,\mu)\,d\mu(x).
$$
We subtract the expression with $f,g$ exchanged. The Hessian terms cancel by symmetry, and the remaining terms yield \eqref{operator-bracket} by \eqref{bracket-def}. All integrals converge by \eqref{field-growth} and \eqref{field-derivative-bounds}.
\end{proof}

The operator identity also determines the bracket among continuous fields. Indeed, if a continuous field $v$ satisfies $\int\nabla \varphi\cdot v(\cdot,\mu)\,d\mu=0$ for every $\varphi\in C_c^\infty$ and every $\mu$, we can prescribe the gradients of $\varphi$ independently at the finitely many atoms of an empirical law. Taking them equal to $v$ at those atoms, we get $v^N=0$, and Lemma~\ref{lem:empirical-faithfulness} applies. The quantifier over all laws is essential here.

\subsection{Observables and the Leibniz rule}
\label{subsec:anchor}

Let $\mathcal{A}_L$ be the algebra of scalar functionals which are finite sums of
$$
F(\mu)=\int_{(\R^d)^m}k(y_1,\ldots,y_m)\,d\mu^{\otimes m}, \qquad \vert \partial_{y_1}^{\beta_1}\cdots\partial_{y_m}^{\beta_m}k\vert \leq C_{\boldsymbol\beta} \prod_{j=1}^m\langle y_j\rangle^{-\vert \beta_j\vert },
$$
where $k$ is smooth and symmetric, and $m=0$ gives the constants. Products are obtained by concatenating and symmetrizing kernels. In particular, $F_\varphi(\mu)=\int \varphi\,d\mu$ belongs to $\mathcal{A}_L$ for $\varphi\in C_c^\infty(\R^d)$.

\begin{proposition}\label{prop:anchor}
For every $f\in\SL$ and $F\in\mathcal{A}_L$, we have $\mathscr{L}_fF\in\mathcal{A}_L$ and $Ff\in\SL$. For $f,g\in\SL$ and $F\in\mathcal{A}_L$,
\begin{align}
[\mathscr{L}_f,\mathscr{L}_g]F&=\mathscr{L}_{\br{f}{g}}F,\label{operator-bracket-AL}\\
\br{f}{Fg}&=F\br{f}{g}+(\mathscr{L}_fF)\,g.\label{anchor-leibniz}
\end{align}
Moreover, if $\mathscr{L}_fF=0$ for every $F\in\mathcal{A}_L$, then $f=0$.
\end{proposition}

\begin{proof}
Multiplying kernels shows that $Ff\in \SL$. Let $F$ and $f$ have kernels $k(\boldsymbol y)$ and $K_f(x;\boldsymbol z)$. Differentiation under the integral represents $\mathscr{L}_fF$ by the symmetrization of $\sum_{i=1}^m \nabla_{y_i}k(\boldsymbol y)\cdot K_f(y_i;\boldsymbol z)$. The weights $\langle y_i\rangle^{-1}$ and $\langle y_i\rangle$ cancel. Further derivatives satisfy the same zero-order estimates, by the argument in \eqref{kernel-weight-cancellation}. Hence this is a kernel in $\mathcal{A}_L$.

The Lions product rule shows that $\mathscr{L}_f$ is a derivation and that $\mathscr{L}_{Ff}=F\mathscr{L}_f$. For $F\in\mathcal{A}_L$, we set $F^N(\boldsymbol x)=F(\mu_{\boldsymbol x}^N)$. Empirical differentiation gives $DF^N f^N=(\mathscr{L}_fF)^N$, as in \eqref{particle-total-derivative}. Applying the classical commutator identity and \eqref{particle-bracket-compatibility}, and then using the density of empirical measures, we obtain \eqref{operator-bracket-AL}. If $\mathscr{L}_fF=0$ for every $F\in\mathcal{A}_L$, then $\mathscr{L}_fF_\varphi=0$ for every $\varphi\in C_c^\infty$, and the argument following Proposition~\ref{prop:operator} implies $f=0$. Finally, substituting $\partial_\mu(Fg)=F\partial_\mu g+g\otimes\partial_\mu F$ in \eqref{bracket-def} proves \eqref{anchor-leibniz}.
\end{proof}

\begin{remark}\label{rem:Lie-Rinehart}
Formula \eqref{anchor-leibniz} is the analog of $[X,\varphi Y]=\varphi[X,Y]+(X\varphi)Y$ for vector fields $X,Y$ and functions $\varphi$ on a manifold. In algebraic terms, $(\mathcal{A}_L,\SL,\rho)$ with $\rho(f)=\mathscr{L}_f$ is a Lie-Rinehart algebra (see \cite[Definition~1.1]{LaurentGengouxLouis2022}), just as smooth functions and vector fields on a manifold, with the canonical action of vector fields as derivations. Lie-Rinehart algebras are the algebraic counterpart of Lie algebroids (see \cite{Meinrenken2024}).
\end{remark}

The class $\SL$ cannot be replaced by the space of all smooth nonlocal vector fields, because $\mathscr{L}_f$ must map $\mathcal{A}_L$ into itself. Let $v,e\in\R^d$ and $\varphi,\psi\in C_c^\infty(\R^d)$. For $f(x,\mu)=\exp\big(\int\varphi\,d\mu\big)v$, the functional $\mathscr{L}_fF_\psi(\mu)=\exp\big(\int\varphi\,d\mu\big)\int\nabla\psi\cdot v\,d\mu$ is in general not polynomial in $\mu$, as can be seen on the laws $(1-s)\delta_a+s\delta_b$, on which the elements of $\mathcal{A}_L$ are polynomial in $s$. For $f(x,\mu)=\big(\int y\cdot e\,d\mu(y)\big)v$, the functional $\mathscr{L}_fF_\psi$ is in general unbounded on $\Ptwo$, whereas the elements of $\mathcal{A}_L$ are bounded. The estimates \eqref{kernel-symbol} are designed so that the weights cancel in the kernel $\nabla k(y)\cdot K_f(y;\boldsymbol z)$ of $\mathscr{L}_fF$, as in the proof of Proposition~\ref{prop:anchor}. In this sense, $\SL$ is a class of symbols with values in nonlocal vector fields, adapted to the algebra $\mathcal{A}_L$ of observables.

In \cite{DelloSchiavoPinamonti2026}, vector fields of the form $V(x,\mu)=\sum_iu_i(\mu)X_i(x)$ on the Wasserstein space over a closed manifold, where the $u_i$ are cylinder functions, that is, smooth functions of finitely many integrals $\int\varphi_j\,d\mu$, and the $X_i$ are smooth vector fields, are bracketed through the commutator of the derivations that they induce on cylinder functions. On $\R^d$, for two such fields $V$ and $W$ satisfying Assumption~\ref{ass:regularity}, Proposition~\ref{prop:operator} and the chain rule show that this commutator is induced by $\br{V}{W}$. The Lie-Lions bracket therefore provides an explicit formula for this bracket, valid for general nonlocal vector fields, and Theorems~\ref{thm:commutation} and~\ref{thm:quantitative} relate it to the flows.

For linear observables, the derivation $F\mapsto\mathscr{L}_fF$ has the distributional representation
$$
\mathscr{L}_fF_\varphi(\mu) =\int\nabla \varphi\cdot f(\cdot,\mu)\,d\mu =\big\langle-\operatorname{div}(f(\cdot,\mu)\mu),\varphi\big\rangle.
$$
At a fixed law $\mu$, the map $F\mapsto\mathscr{L}_fF(\mu)$ only depends on the velocity that $f$ imposes on $\mu$. The next proposition identifies this velocity with a tangent vector in the sense of Otto's calculus, which relates the Lions calculus used in this paper to the Riemannian formalism on $\Ptwo$ (see \cite{AGS,Otto2001}).

Let $\mu\in\Ptwo$. We set
$$
T_\mu\mathcal{P}_2 =\overline{\{\nabla \varphi\ \mid\ \varphi\in C_c^\infty(\R^d)\}}^{L^2(\mu)}, \qquad P_\mu:L^2(\mu;\R^d)\longrightarrow T_\mu\mathcal{P}_2.
$$
This is the Otto tangent space, and $P_\mu$ denotes its orthogonal projection (see \cite[Chapter~8]{AGS}).

\begin{proposition}
\label{prop:otto-representative}
For every $f\in\SL$ and every $\mu\in\Ptwo$, we have
\begin{equation*}
\mathscr{L}_fF(\mu) =\int\partial_\mu F(\mu)(x)\cdot P_\mu f(\cdot,\mu)(x)\,d\mu(x) \qquad\forall F\in\mathcal{A}_L.
\end{equation*}
\end{proposition}

\begin{proof}
The derivative $\partial_\mu F(\mu)$ is the gradient of a smooth bounded function with bounded gradient. Multiplication of its potential by compactly supported cutoffs shows that it belongs to $T_\mu\mathcal{P}_2$, and the formula follows by orthogonality.
\end{proof}

By orthogonality, $\operatorname{div}(\mu P_\mu f(\cdot,\mu))=\operatorname{div}(\mu f(\cdot,\mu))$: the fields $f(\cdot,\mu)$ and $P_\mu f(\cdot,\mu)$ move $\mu$ in the same way at first order, and $P_\mu f(\cdot,\mu)$ is the velocity of minimal norm with this property. Hence $\Vert P_{\mu_t}f(\cdot,\mu_t)\Vert_{L^2(\mu_t)}$ is the metric speed in $W_2$ of the curve $\mu_t=\Phi_t^f\mu$, for almost every $t$ (see \cite[Theorem~8.3.1 and Proposition~8.4.5]{AGS}). For instance, in dimension two, the rotation field $f(x)=Jx$, where $J$ is the rotation by $\pi/2$, moves every point except the origin, but leaves invariant every rotation-invariant law $\mu$, and $P_\mu f=0$ for such laws. This explains why the last assertion of Proposition~\ref{prop:anchor} requires all laws: at a fixed law, $P_\mu f(\cdot,\mu)$ may vanish even when $f(\cdot,\mu)$ does not.

\section{Applications to population dynamics}\label{sec:examples}

In the following examples, we compare the Lie-Lions bracket with the frozen spatial bracket $[f,g]_x$. In Section~\ref{subsec:symmetries}, the flows commute although the frozen bracket does not vanish. In Sections~\ref{subsec:skater}, \ref{subsec:flock-turning} and~\ref{subsec:congestion}, the frozen bracket vanishes for every law, whereas the Lie-Lions bracket does not. Both brackets are nonzero in Sections~\ref{subsec:kuramoto} and~\ref{subsec:three-body}, where only the Lie-Lions bracket depends on the variation of the coherence and on the relative velocities of the interacting agents, respectively.

\subsection{Symmetries}\label{subsec:symmetries}

Most interaction fields of collective dynamics depend only on the relative positions of the agents, and are thus invariant under translations of the whole population. We characterize such symmetries by the vanishing of the Lie-Lions bracket with the generator of the motion, as for vector fields on a manifold, where a field is invariant under the flow of another field if and only if their bracket vanishes. The frozen spatial bracket does not vanish for these fields in general, so that its vanishing is not necessary for the commutation of the flows. This characterization is used again in Sections~\ref{subsec:kuramoto} and~\ref{subsec:three-body}.

Let $A\in\R^{d\times d}$ and $v\in\R^d$. The local field $g(x)=Ax+v$ belongs to $\SL$, and its flow acts on laws by $\Phi_s^g\mu=(T_s)_\#\mu$, where $T_s(x)=e^{sA}x+\int_0^se^{rA}v\,dr$. When $A$ is skew-symmetric, $T_s$ is a rigid motion of the whole population.

\begin{proposition}\label{prop:symmetries}
Let $f$ satisfy Assumption~\ref{ass:regularity}, and let $g(x)=Ax+v$. Then $\br{f}{g}=0$ if and only if
$$
f(T_s(x),(T_s)_\#\mu)=e^{sA}f(x,\mu) \qquad\forall s\in\R \qquad\forall x\in\R^d \qquad\forall\mu\in\Ptwo.
$$
\end{proposition}

\begin{proof}
We have $\mathcal{D}_fg=Af$, since $\partial_xg=A$ and $\partial_\mu g=0$. Along $x_s=T_s(x)$ and $\mu_s=(T_s)_\#\mu$, the chain rule recalled after \eqref{total-derivative} gives $\frac{d}{ds}f(x_s,\mu_s)=\mathcal{D}_gf(x_s,\mu_s)$. Hence
$$
\frac{d}{ds}\big(e^{-sA}f(x_s,\mu_s)\big)=-e^{-sA}\br{f}{g}(x_s,\mu_s),
$$
and the equivalence follows.
\end{proof}

The frozen spatial bracket is $[f,g]_x(x,\mu)=Af(x,\mu)-\partial_xf(x,\mu)(Ax+v)$. It differs from $\br{f}{g}(x,\mu)$ by $\int\partial_\mu f(x,\mu)(y)(Ay+v)\,d\mu(y)$, and it does not vanish for equivariant fields in general. For instance, let $A=0$. A field depending only on relative positions, i.e., such that $f(x+v,(\tau_v)_\#\mu)=f(x,\mu)$ for every $v\in\R^d$, with $\tau_v(y)=y+v$, satisfies $\br{f}{v}=0$, and thus its flow commutes with every translation of the population. This is the case of the convolution fields $f(x,\mu)=\int\psi(x-y)\,d\mu(y)$ with $\psi\in C_b^\infty(\R^d;\R^d)$, and of the aggregation field of Section~\ref{subsec:skater}, whereas $[f,v]_x(x,\mu)=-\partial_xf(x,\mu)v$. Similarly, in dimension two, let $g(x)=Jx$, where $J$ is the rotation of angle $\pi/2$. If $\psi\in C_b^\infty(\R^2;\R^2)$ satisfies $\psi(Rz)=R\psi(z)$ for every rotation $R$, then the corresponding convolution field satisfies $\br{f}{g}=0$, although $[f,g]_x(x,\mu)=Jf(x,\mu)-\partial_xf(x,\mu)Jx$ does not vanish in general. Equivariant nonlocal fields are also considered in \cite[Section~6]{GoncharovaPogodaevStaritsyn2026}, in relation with projected controllability.

\subsection{Aggregation and collective rotation}
\label{subsec:skater}

\begin{figure}[H]
\centering
\begin{tikzpicture}[x=1cm,y=1cm,>=Stealth,font=\small]
  \definecolor{popblue}{RGB}{31,67,112}
  \definecolor{popred}{RGB}{156,47,50}
  \newcommand{\skater}[3]{
    \begin{scope}[shift={(#1,4.75)}]
      \draw[black!45,dashed,line width=.45pt] (0,-.2)--(0,3.75);
      \begin{scope}[draw=#2,line width=1.15pt,line cap=round,line join=round]
        \draw[fill=white] (0,2.83) circle (.2);
        \draw (0,2.62)--(0,1.45);
        \draw (0,1.45)--(.04,.14);
        \draw[line width=.9pt] (-.22,.06)--(.32,.06);
        \draw (0,1.45)--(.43,.86)--(.07,.72);
        \ifnum#3=0
          \draw (0,2.38)--(-.62,2.33)--(-1.12,2.27);
          \draw (0,2.38)--(.62,2.33)--(1.12,2.27);
        \else
          \draw (0,2.38)--(-.42,2.95)--(-.07,3.55);
          \draw (0,2.38)--(.42,2.95)--(.07,3.55);
        \fi
      \end{scope}
    \end{scope}}
  \skater{2.0}{popblue}{0}
  \skater{8.0}{popred}{1}
  \draw[->,popblue,line width=.9pt]
    (2,6.65)++(225:.85 and .24) arc[start angle=225,end angle=310,x radius=.85,y radius=.24];
  \draw[->,popred,line width=.9pt]
    (8,6.65)++(150:.85 and .24) arc[start angle=150,end angle=410,x radius=.85,y radius=.24];
  \node[align=center,font=\scriptsize] at (2.0,4.3) {Arms extended};
  \node[align=center,font=\scriptsize] at (8.0,4.3) {Arms raised above the head};
    \newcommand{\skaterpopulation}[4]{%
    \begin{scope}[shift={(#1,1.8)}]
      \begin{scope}[scale=#2]
        \foreach \x/\y in {1.2/.05,.35/.75,-.9/.4,-.8/-.6,.15/-.6}
          \node[circle,draw=#3,fill=#4,line width=.7pt,inner sep=1.75pt]
            at (\x,\y) {};
      \end{scope}
      \draw[black!65,line width=.45pt]
        (-.075,0)--(.075,0) (0,-.075)--(0,.075);
    \end{scope}}
  \node at (2.0,3.7) {Spread population};
  \node at (8.0,3.7) {Concentrated population};
  \skaterpopulation{2.0}{1}{popblue}{white}
  \skaterpopulation{8.0}{.34}{popred}{popred}
  \draw[->,popblue,line width=.9pt]
    (2,1.8)++(20:1.45) arc[start angle=20,end angle=85,radius=1.45];
  \draw[->,popred,line width=.9pt]
    (8,1.8)++(-40:1.45) arc[start angle=-40,end angle=180,radius=1.45];
  \node[align=center,font=\scriptsize] at (2.0,.3)
    {Larger inertia\\slower rotation};
  \node[align=center,font=\scriptsize] at (8.0,.3)
    {Smaller inertia\\faster rotation};
\end{tikzpicture}
\caption{Top: a skater spins faster with the arms raised above the head, since the moment of inertia about the rotation axis decreases while the angular momentum is conserved. Bottom: for fixed $L>0$, contraction reduces the inertia $V$ and increases the angular speed $a(V)=L/(1+V)$. The curved arrows indicate the relative speeds schematically.}
\label{fig:skater}
\end{figure}
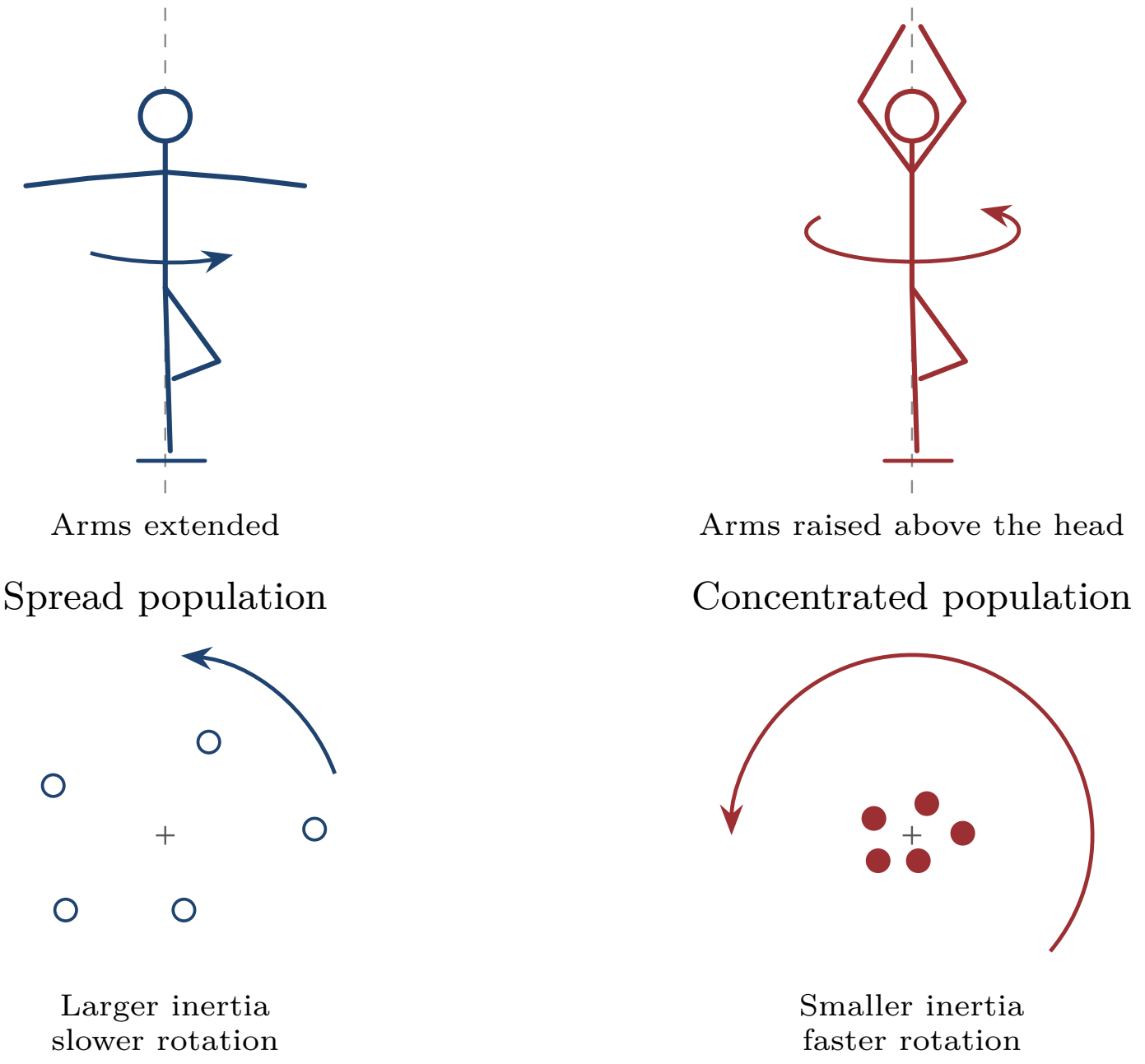

We consider the aggregation velocity $f=-\nabla W*\mu$ associated with the quadratic potential $W(z)=k\vert z\vert^2/2$ \cite{CarrilloDiFrancescoFigalliLaurentSlepcev2011}. It is directed towards the barycenter $m(\mu)=\int x\,d\mu(x)$. We combine it with rotation at an angular speed depending on the variance, as suggested by the increase in a skater's angular speed when the mass moves towards the rotation axis \cite{Perez2016} (see Figure~\ref{fig:skater}).

For a planar mass distribution, the variance is the moment of inertia per unit mass about the barycenter. Adding a fixed positive contribution, normalized to one, we obtain the regularized angular-speed law $a(V)=L/(1+V)$, with prescribed $L\in\R$.

In dimension two, we set
$$
J=\begin{pmatrix}0&-1\\1&0\end{pmatrix},\qquad R_\theta=e^{\theta J},\qquad V(\mu)=\int\vert x-m(\mu)\vert^2\,d\mu(x),
$$
and, for $k>0$, we define
$$
f(x,\mu)=-k(x-m(\mu)),\qquad g(x,\mu)=a(V(\mu))J(x-m(\mu)),\qquad a(v)=\frac{L}{1+v}.
$$
The matrices $-kI$ and $a(V)J$ commute, and thus $[f,g]_x=0$ for every fixed law. The maps $\widetilde f$ and $\widetilde g$ are smooth, and the identities $\partial_\mu m(\mu)(y)=I$ and $\partial_\mu V(\mu)(y)=2(y-m(\mu))$ imply the derivative bounds of Assumption~\ref{ass:regularity} on second-moment balls. The field $f$ does not belong to $\SL$: by \eqref{symbol-first-bounds}, at fixed $x$, every field of $\SL$ is bounded with respect to the law, whereas $f(x,\mu)$ is unbounded as $m(\mu)$ varies. The same holds for $g$ when $L\neq0$.

Both flows preserve $m=m(\mu)$. Contraction changes the variance to $e^{-2ks}V$, whereas rotation preserves it. Writing $V=V(\mu)$,
$$
\Phi_s^f\mu=(x\mapsto m+e^{-ks}(x-m))_\#\mu,\qquad \Phi_t^g\mu=(x\mapsto m+R_{t a(V)}(x-m))_\#\mu.
$$
Since $\mathscr{L}_fV=-2kV$ and $\mathscr{L}_gV=0$, formula \eqref{bracket-def} gives
\begin{equation}\label{skater-bracket}
\br{f}{g}(x,\mu) =-2kV a'(V)J(x-m) =\frac{2kLV}{(1+V)^2}J(x-m).
\end{equation}
The two compositions $\Phi_t^g\Phi_s^f$ and $\Phi_s^f\Phi_t^g$ have the same final variance, but their angles differ by $\theta(s,t)=t(a(e^{-2ks}V)-a(V))$. For $L>0$, $V>0$ and $s,t>0$, aggregation followed by rotation gives a larger angle, since the rotation takes place at smaller inertia. The flow commutator is
$$
\Phi_{-t}^g\circ\Phi_{-s}^f\circ \Phi_t^g\circ\Phi_s^f\mu =(x\mapsto m+R_{\theta(s,t)}(x-m))_\#\mu.
$$
We recover \eqref{skater-bracket} from $\partial_s\theta(0,t)=-2ktVa'(V)$. The commutator changes the law whenever it is not invariant under this rotation. The flow of $\br{f}{g}$ is the rotation about $m$ at the angular speed $-2kVa'(V)$, which preserves $m$ and $V$. The $W_2$ distance between the images of $\mu$ under two rotations about $m$, of angles $\alpha$ and $\beta$, is at most $\vert\alpha-\beta\vert V^{1/2}$, and a Taylor expansion of $a$ gives $\vert\theta(t,t)+2kVa'(V)t^2\vert\leq C\vert t\vert^3$ for $\vert t\vert\leq1$, with $C$ depending only on $k$ and $L$. Hence $W_2\big(\mathcal{C}_t^{f,g}\mu,\exp(t^2\br{f}{g})\mu\big)\leq C\vert t\vert^3V^{1/2}$, and the first estimate of Theorem~\ref{thm:quantitative} holds for these fields, although $f\notin\SL$.

\subsection{Synchronization and coherence-dependent frequency}\label{subsec:kuramoto}

In models of synchronization, the coherence of a population of oscillators can play the role of the variance in Section~\ref{subsec:skater}. The phase of an oscillator is represented by $\theta\in\R$. We consider the coupling of the Kuramoto model \cite{AcebronBonillaPerezVicenteRitortSpigler2005} and a common rotation whose frequency depends on the coherence:
$$
f(\theta,\mu)=K\int\sin(\phi-\theta)\,d\mu(\phi),\qquad g(\theta,\mu)=\omega_0+\omega_1R(\mu)^2,\qquad R(\mu)^2=\iint\cos(\phi-\phi')\,d\mu(\phi)\,d\mu(\phi'),
$$
with $K>0$ and $\omega_0,\omega_1\in\R$. Here $R(\mu)=\big\vert\int e^{i\phi}\,d\mu(\phi)\big\vert\in[0,1]$ is the order parameter. Both fields belong to $\mathcal{B}_L$.

Since $g$ does not depend on $\theta$ and $R$ is invariant under common phase shifts, the flow of $g$ shifts all phases by $tg(\mu)$. The coupling $f$ is invariant under such shifts, and thus $\mathcal{D}_gf=0$, as in the proof of Proposition~\ref{prop:symmetries}. Moreover, $\partial_\mu(R^2)(\mu)(\phi)=2f(\phi,\mu)/K$. We infer that $\mathcal{D}_fg=\omega_1\mathscr{L}_f(R^2)$ and
\begin{equation}\label{kuramoto-bracket}
\br{f}{g}(\theta,\mu)=\frac{2\omega_1}{K}\int f(\phi,\mu)^2\,d\mu(\phi).
\end{equation}
The bracket \eqref{kuramoto-bracket} comes only from the nonlocal term $\omega_1\mathscr{L}_f(R^2)$, that is, from the effect of the coupling on the coherence that sets the frequency of $g$. At a fixed law with $g(\mu)\neq0$, the bracket is not a new direction, since it is a common phase velocity like $g$. Its amplitude $\omega_1\frac{d}{ds}R(\Phi_s^f\mu)^2\vert_{s=0}$ is the rate at which synchronization accelerates the common rotation. In contrast, the frozen spatial bracket is $[f,g]_x(\theta,\mu)=Kg(\mu)\int\cos(\phi-\theta)\,d\mu(\phi)$. It does not vanish when $g(\mu)\neq0$ and $R(\mu)>0$, even for $\omega_1=0$, although the flows then commute by Theorem~\ref{thm:commutation}.

Since common phase shifts preserve $R$ and commute with the flow of $f$, the order in which coupling and frequency modulation are applied does not affect the coherence. The two compositions $\Phi_t^g\Phi_s^f\mu$ and $\Phi_s^f\Phi_t^g\mu$ differ by a common phase shift of angle $t\omega_1\big(R(\Phi_s^f\mu)^2-R(\mu)^2\big)$, which is nonnegative for $\omega_1,s,t\geq0$, since $\mathscr{L}_f(R^2)=\frac{2}{K}\int f^2\,d\mu\geq0$. Synchronizing before rotating therefore increases the frequency of the rotation, and the derivative of this angle at $s=0$ is $t\br{f}{g}(\theta,\mu)$. For the same reason, $\Phi_t^{f+g}\mu$ is obtained from $\Phi_t^f\mu$ by a common phase shift. Hence the Lie splitting $\exp(tg)\circ\exp(tf)$ reproduces the evolution of $R$ exactly, and its error is a common phase shift, whose leading term is given by $\frac{t}{2}\br{f}{g}$ (see Remark~\ref{rem:splitting}).

\subsection{Transmission of a turn in a flock}\label{subsec:flock-turning}
Leader-follower hierarchies have been inferred from delays between directional changes in pigeon flocks \cite{NagyAkosBiroVicsek2010}. We consider the simplest such hierarchy: a leader, a relay that aligns with the leader, and a follower that aligns with the relay. The two alignment interactions act on different groups, so that their frozen spatial fields commute for every law, and the classical Lie bracket predicts that the order of the interactions is irrelevant. We show that their Lie-Lions bracket is nonzero, and that it describes the transmission of a turn from the relay to the follower: it explains why the order matters, and it gives the delay, measured by the order in time, with which a turn of the leader reaches the follower. Our simplified model combines directed interactions \cite{MotschTadmor2011} and alignment towards an observed heading \cite{DegondMotsch2008}, with a follower group, a relay group and a leader group.

The state $x=(s,\theta)\in\R^2$ consists of a fixed group label $s$ and a heading angle $\theta\in\R$. Labels $1,2,3$ correspond respectively to the follower, relay and leader. The interaction coefficients are periodic in the angle, which is represented on the real line. We describe only changes of heading, with no spatial motion.

Let $\chi_i\in C_c^\infty(\R;[0,1])$ have pairwise disjoint supports and satisfy $\chi_i(j)=\delta_{ij}$ for $i,j\in\{1,2,3\}$. For $\kappa>0$, we set
\begin{equation*}
\begin{aligned}
f((s,\theta),\mu) &=\Big(0,\, 3\kappa\chi_2(s) \int\chi_3(r)\sin(\phi-\theta)\,d\mu(r,\phi)\Big),\\
g((s,\theta),\mu) &=\Big(0,\, 3\kappa\chi_1(s) \int\chi_2(r)\sin(\phi-\theta)\,d\mu(r,\phi)\Big).
\end{aligned}
\end{equation*}
The sine is the tangential component of the observed direction at the current heading. The labels remain fixed. When each label has mass $1/3$, the factor $3$ normalizes the observed group mass. The field $f$ turns the relay towards the leader, leaving the leader and follower fixed, while $g$ turns the follower towards the fixed relay. These fields are instances of the multiple-wise Kuramoto model with distinguishable agents of \cite[Example~2.7]{PaulRossiTrelat2025}, with pairwise interactions ($m=1$), zero natural frequencies and hierarchical weights $3\kappa\chi_2(s)\chi_3(r)$ and $3\kappa\chi_1(s)\chi_2(r)$.

The kernels and all their derivatives are bounded, and thus Assumption~\ref{ass:regularity} holds. Since the fields have only an angular component and act on disjoint label supports, $[f,g]_x=0$ for every law. The Lie-Lions bracket has zero label component and angular component
$$
\br{f}{g}_\theta((s,\theta),\mu) =9\kappa^2\chi_1(s)\int\chi_2(r)^2\cos(\phi-\theta) \int\chi_3(u)\sin(\psi-\phi)\,d\mu(u,\psi) \,d\mu(r,\phi).
$$
Indeed, only $\mathscr{L}_fg$ remains in \eqref{bracket-def}: $g$ depends on the heading of the relay through the law, and $f$ changes this heading. One factor $\chi_2(r)$ is the relay weight in $g$, and the other comes from the action of $f$ on the relay.

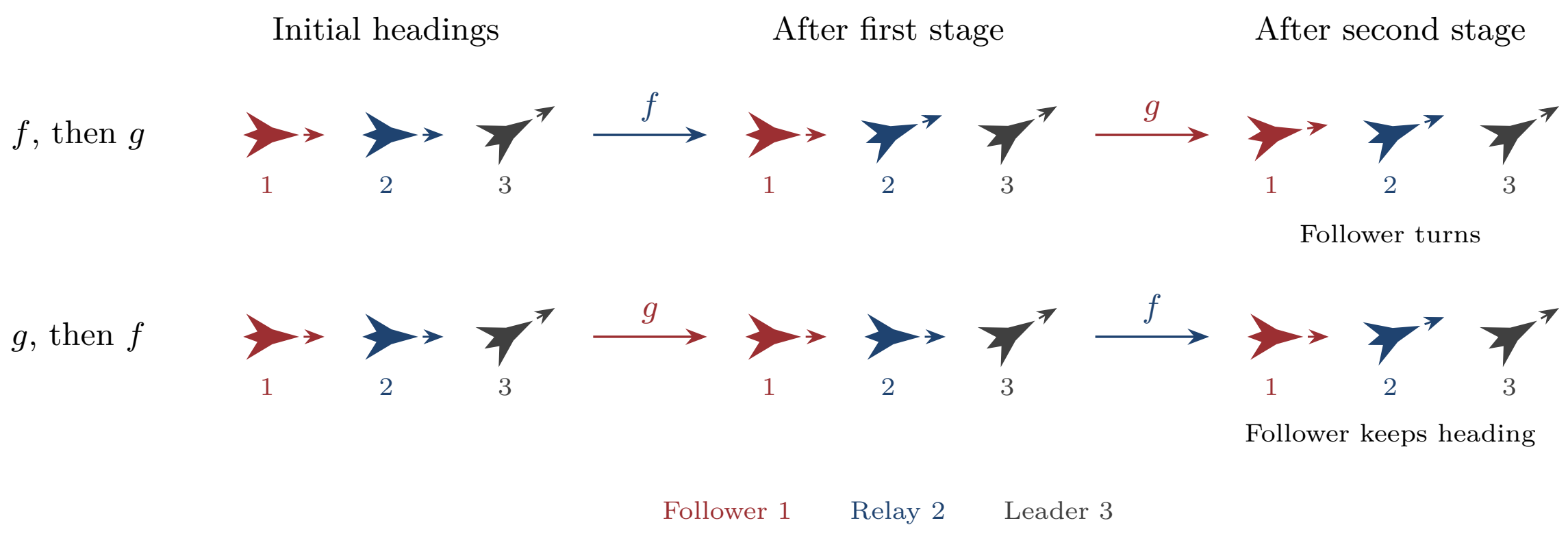
\begin{figure}[H]
\centering
\begin{tikzpicture}[x=1cm,y=1cm,>=Stealth,font=\small]
  \definecolor{popblue}{RGB}{31,67,112}
  \definecolor{popred}{RGB}{156,47,50}
  \newcommand{\flockbird}[3]{%
    \begin{scope}[shift={(#1)},rotate=#2,scale=.72]
      \fill[#3] (.43,0)--(.07,.10)--(-.28,.31)--(-.13,.07)
        --(-.32,0)--(-.13,-.07)--(-.28,-.31)--(.07,-.10)--cycle;
      \draw[->,#3,line width=.55pt] (.49,0)--(.77,0);
    \end{scope}}
  \newcommand{\flockstage}[4]{%
    \begin{scope}[shift={(#1)}]
      \flockbird{-1.15,0}{#2}{popred}
      \flockbird{0,0}{#3}{popblue}
      \flockbird{1.15,0}{#4}{black!75}
      \node[font=\scriptsize,text=popred] at (-1.15,-.48) {$1$};
      \node[font=\scriptsize,text=popblue] at (0,-.48) {$2$};
      \node[font=\scriptsize,text=black!75] at (1.15,-.48) {$3$};
    \end{scope}}
  \node at (2.15,3.6) {Initial headings};
  \node at (7,3.6) {After first stage};
  \node at (11.85,3.6) {After second stage};
  \node[anchor=east] at (-.05,2.6) {$f$, then $g$};
  \node[anchor=east] at (-.05,.65) {$g$, then $f$};
  \flockstage{2.15,2.6}{0}{0}{30}
  \flockstage{7,2.6}{0}{20}{30}
  \flockstage{11.85,2.6}{10}{20}{30}
  \flockstage{2.15,.65}{0}{0}{30}
  \flockstage{7,.65}{0}{0}{30}
  \flockstage{11.85,.65}{0}{20}{30}
  \draw[->,popblue,line width=.65pt] (4.15,2.6)--node[above] {$f$}(5.25,2.6);
  \draw[->,popred,line width=.65pt] (9,2.6)--node[above] {$g$}(10.1,2.6);
  \draw[->,popred,line width=.65pt] (4.15,.65)--node[above] {$g$}(5.25,.65);
  \draw[->,popblue,line width=.65pt] (9,.65)--node[above] {$f$}(10.1,.65);
  \node[font=\scriptsize] at (11.85,1.65) {Follower turns};
  \node[font=\scriptsize] at (11.85,-.3) {Follower keeps heading};
  \node[font=\scriptsize] at (7,-1.05)
    {\textcolor{popred}{Follower $1$}\qquad
     \textcolor{popblue}{Relay $2$}\qquad
     \textcolor{black!75}{Leader $3$}};
\end{tikzpicture}
\caption{Applying $f$ before $g$ turns the relay towards the leader, then the follower towards the relay. In the reverse order, $g$ leaves the initially aligned follower unchanged, and $f$ turns only the relay.}
\label{fig:flock}
\end{figure}

When each group has a single heading, we set $\mu_{\boldsymbol\theta} =\frac{1}{3}\sum_{i=1}^3\delta_{(i,\theta_i)}$. The flow of $f+g$ preserves this family and is given by
$$
\dot\theta_1=\kappa\sin(\theta_2-\theta_1),\qquad \dot\theta_2=\kappa\sin(\theta_3-\theta_2),\qquad \dot\theta_3=0.
$$
For these atomic laws, the bracket is
$$
\br{f}{g}((s,\theta),\mu_{\boldsymbol\theta}) =\big(0,\kappa^2\chi_1(s)\cos(\theta_2-\theta) \sin(\theta_3-\theta_2)\big).
$$
Let us assume that $\theta_1(0)=\theta_2(0)=\theta_0$ and $\theta_3(0)=\theta_0+\delta$, that is, the follower and the relay are initially aligned and the leader has turned. Applying $f$ before $g$ gives follower heading $\theta_0+\kappa^2\sin\delta\,t^2+\mathrm{O}(t^3)$, while the reverse order leaves it exactly at $\theta_0$ (see Figure~\ref{fig:flock}). The corrections $\pm\frac{t}{2}\br{f}{g}$ respectively add and cancel the leading turn generated by $f+g$.

Along the flow of $f+g$, differentiating the follower's turning rate, we get
\begin{equation}\label{flock-angular-acceleration}
\ddot\theta_1 =\underbrace{\kappa\cos(\theta_2-\theta_1)\dot\theta_2} _{\text{relay's turn: }\br{f}{g}_\theta} -\underbrace{\kappa\cos(\theta_2-\theta_1)\dot\theta_1} _{\text{follower's own turn}}.
\end{equation}
The first term in \eqref{flock-angular-acceleration} is the angular component of $\br{f}{g}((1,\theta_1),\mu_{\boldsymbol\theta})$, namely the contribution of the relay's change of heading to the follower's acceleration. For the initial data above, we have
$$
\dot\theta_1(0)=0,\qquad \ddot\theta_1(0)=\kappa^2\sin\delta,\qquad \theta_1(t)=\theta_0+\frac{\kappa^2}{2}\sin\delta\,t^2+\mathrm{O}(t^3).
$$
The follower thus starts to turn at the second order in time, with an acceleration given by the bracket. Since the kernels of $f,g$ belong to $C_b^\infty$, Proposition~\ref{prop:bounded-kernels} gives constants $C,t_0>0$ such that, for every $\mu \in \mathcal{P}_2(\R^2)$ and $\vert t\vert \leq t_0$,
\begin{align*}
W_2\big(\Phi_t^g\Phi_t^f\mu,\Phi_t^{f+g+\frac{t}{2}\br{f}{g}}\mu\big) &\leq C\vert t\vert^3,\\
W_2\big(\Phi_t^f\Phi_t^g\mu,\Phi_t^{f+g-\frac{t}{2}\br{f}{g}}\mu\big) &\leq C\vert t\vert^3.
\end{align*}
These estimates are uniform in the initial law and the population size, although the periodic kernels do not satisfy the decay assumptions defining $\SL$.

The same computation applies to hierarchies of arbitrary depth. Let $n\geq2$, and let $\chi_1,\ldots,\chi_{n+1}\in C_c^\infty(\R;[0,1])$ have pairwise disjoint supports and satisfy $\chi_i(j)=\delta_{ij}$. For $1\leq k\leq n$, we set
$$
f_k((s,\theta),\mu)=\Big(0,\,(n+1)\kappa\chi_k(s)\int\chi_{k+1}(r)\sin(\phi-\theta)\,d\mu(r,\phi)\Big).
$$
The field $f_k$ turns the group $k$ towards the group $k+1$, and the group $n+1$ is the leader. For $n=2$, we have $f_1=g$ and $f_2=f$. As above, the frozen spatial brackets of these fields vanish. We now set $\mu_{\boldsymbol\theta}=\frac{1}{n+1}\sum_{i=1}^{n+1}\delta_{(i,\theta_i)}$, together with $B_1=f_1$ and $B_k=\br{f_k}{B_{k-1}}$ for $2\leq k\leq n$. By Proposition~\ref{prop:particle}, these brackets can be computed in the particle realization with one particle per group, in which $f_k$ becomes the field $\kappa\sin(\theta_{k+1}-\theta_k)$ in the $k$th heading. An induction on $k$ gives
$$
B_k((1,\theta_1),\mu_{\boldsymbol\theta})=\Big(0,\,\kappa^k\prod_{j=1}^{k-1}\cos(\theta_{j+1}-\theta_j)\,\sin(\theta_{k+1}-\theta_k)\Big).
$$
Let us assume that the groups $1,\ldots,n$ have the same initial heading $\theta_0$, and that the leader has the heading $\theta_0+\delta$. Then $B_k$ vanishes at the follower for $k<n$, and the flow of $f_1+\cdots+f_n$ gives
$$
\theta_1(t)=\theta_0+\frac{\kappa^n}{n!}\sin\delta\,t^n+\mathrm{O}(t^{n+1}),
$$
where $\kappa^n\sin\delta$ is the angular component of $B_n$ at the follower. For hierarchies of depth $n$, all the frozen spatial fields still commute, whereas the Lie-Lions brackets of length $k\leq n$ transmit the turn of the leader down $k$ levels. The turn reaches the group at distance $n$ at order $n$ in time, through a bracket of length $n$. This is the nonlocal counterpart of a classical fact in geometric control, where a direction obtained from brackets of length $k$ is reached at order $k$ in time (see \cite{AgrachevSachkov2004}).

\subsection{Effective three-body interactions}\label{subsec:three-body}

Let us consider two pairwise interactions given by convolution fields, as in aggregation models \cite{CarrilloDiFrancescoFigalliLaurentSlepcev2011} (for instance, an attraction and a repulsion with different ranges):
$$
f(x,\mu)=\int\psi_f(x-y)\,d\mu(y),\qquad g(x,\mu)=\int\psi_g(x-y)\,d\mu(y),
$$
with $\psi_f,\psi_g\in C_b^\infty(\R^d;\R^d)$, so that $f,g\in\mathcal{B}_L$. Since $\partial_\mu g(x,\mu)(z)=-D\psi_g(x-z)$, we have $\mathcal{D}_fg(x,\mu)=\int D\psi_g(x-z)\big(f(x,\mu)-f(z,\mu)\big)\,d\mu(z)$, and therefore
\begin{equation}\label{three-body-bracket}
\br{f}{g}(x,\mu)=\iint\Big(D\psi_g(x-z)\big(\psi_f(x-y)-\psi_f(z-y)\big)-D\psi_f(x-z)\big(\psi_g(x-y)-\psi_g(z-y)\big)\Big)\,d\mu(y)\,d\mu(z).
\end{equation}
This field has measure degree two, i.e., it is a three-body interaction. In the first term, the pair interaction $\psi_g(x-z)$ is differentiated along the relative velocity $\psi_f(x-y)-\psi_f(z-y)$ that a third agent, located at $y$, imposes on the pair through $f$. The frozen spatial bracket involves instead the velocity $f(x,\mu)$ of the agent located at $x$. In particular, it does not vanish when $f$ is replaced with a constant field, contrary to $\br{f}{g}$ (see Section~\ref{subsec:symmetries}). By the Jacobi identity in $\mathcal{B}_L$ and Proposition~\ref{prop:symmetries}, the bracket of two translation-invariant fields of $\mathcal{B}_L$ is translation invariant, as can be read on \eqref{three-body-bracket}.

\medskip

Let us assume that the two interactions alternate on time intervals of length $t$. The agent located at $x$ then reacts through $g$ to agents whose positions have just been changed through $f$, so that its velocity depends on pairs of other agents. By Proposition~\ref{prop:bounded-kernels} and Remark~\ref{rem:splitting}, a step $\exp(tg)\circ\exp(tf)$ of this alternation coincides with the flow of $f+g+\frac{t}{2}\br{f}{g}$ during the time $t$, up to a remainder of order three in $W_2$, uniformly in the number of agents. Compared with the simultaneous action of $f$ and $g$, the only new term at this order is $\frac{t}{2}\br{f}{g}$, which is a three-body interaction by \eqref{three-body-bracket}: the agent at $x$ responds to the agent at $z$ through $g$, along the relative velocity that a third agent at $y$ imposes on the pair through $f$. In particle coordinates, the corresponding velocity of the $i$th agent is a double sum over pairs of other agents. The alternation of two pairwise interactions thus produces an effective three-body interaction, whose amplitude is proportional to the switching period $t$, and which vanishes in the limit of fast switching.

In the terminology of \cite{PaulRossiTrelat2025}, Proposition~\ref{prop:symbol} shows that the bracket of an $m$-wise and an $n$-wise interaction is at most $(m+n)$-wise. In the expansion of an alternation of pairwise interactions, the term of order $t^k$ involves brackets of length $k$, which are at most $k$-wise interactions. Alternating pairwise interactions thus generates multiple-wise interactions of all orders, such as those studied in \cite{PaulRossiTrelat2025}. For instance, the bracket of Section~\ref{subsec:flock-turning} is a $2$-wise interaction between the follower, the relay and the leader.

\subsection{Occupancy-dependent circulation}
\label{subsec:congestion}

In pedestrian models, the velocity of an individual depends on the surrounding population (see \cite{BellomoGibelli2015,CarrilloMartinWolfram2016,ColGarLecPedestrian,PicRosNonlocalVelocity}). We consider a crowd near a point of interest $x_*$, such as a stage or a monument, with two behaviors whose intensities depend on the occupancy perceived near $x_*$: gathering towards $x_*$ and circulation about $x_*$. This is the situation of Section~\ref{subsec:skater}, with the variance replaced by a local occupancy measured by a sensor. Unlike the fields of Section~\ref{subsec:skater}, the fields below belong to $\SL$, so that Theorem~\ref{thm:quantitative} applies, with constants independent of the number of pedestrians.

Let $x_*\in\R^2$ and $\sigma>0$. We set
$$
\eta(y)=\exp\Big(-\frac{\vert y-x_*\vert^2}{2\sigma^2}\Big),\qquad M(\mu)=\int\eta(y)\,d\mu(y).
$$
Since $\mu$ is a probability measure, $M\in(0,1]$ is a weighted occupancy fraction. Let $\kappa(M)=\kappa_0+\kappa_1M$ and $\omega(M)=\omega_0+\omega_1M$, with $\kappa_0>0$, $\kappa_1,\omega_1\geq0$ and $\omega_0\in\R$. We define
$$
f(x,\mu)=-\kappa(M(\mu))(x-x_*),\qquad g(x,\mu)=\omega(M(\mu))J(x-x_*).
$$
The field $f$ describes gathering, with $\kappa_1\geq0$ allowing stronger attraction at larger occupancy. For a fixed law, the matrices $-\kappa(M)I$ and $\omega(M)J$ commute, and thus $[f,g]_x=0$.

Both fields belong to $\SL$. Indeed, $f(x,\mu)=-\kappa_0(x-x_*) -\int\kappa_1(x-x_*)\eta(y)\,d\mu(y)$, and $g$ has an analogous representation. These kernels are affine in $x$, with constant or Gaussian dependence on the measure variable. The Gaussian and its derivatives decay rapidly, so these kernels satisfy the symbol estimates at every order and Proposition~\ref{prop:symbol} applies.

\begin{figure}[H]
\centering
\begin{tikzpicture}[x=1cm,y=1cm,>=Stealth,font=\small]
  \definecolor{popblue}{RGB}{31,67,112}
  \definecolor{popred}{RGB}{156,47,50}
  \node at (1.9,3.55) {Occupancy sensor};
  \begin{scope}[shift={(1.9,1.6)}]
    \fill[gray!5] (0,0) circle[radius=1.4];
    \fill[gray!8] (0,0) circle[radius=1.0];
    \fill[gray!12] (0,0) circle[radius=.6];
    \draw[black!65,line width=.45pt] (-.07,0)--(.07,0) (0,-.07)--(0,.07);
    \node[font=\scriptsize,anchor=north east] at (-.08,-.08) {$x_*$};
    \foreach \ang/\rad in {10/1.3,65/1.32,145/1.16,205/1.3,290/1.32}{%
      \draw[popblue,fill=white,line width=.65pt] (\ang:\rad) circle[radius=2pt];
      \fill[popred] (\ang:{.62*\rad}) circle[radius=2pt];
      \draw[->,black!60,line width=.45pt]
        (\ang:{.9*\rad})--(\ang:{.7*\rad});
    }
  \end{scope}
  \draw[popblue,fill=white,line width=.65pt] (.7,-.03) circle[radius=2pt];
  \node[font=\scriptsize,anchor=west] at (.9,-.03) {initial};
  \fill[popred] (2.2,-.03) circle[radius=2pt];
  \node[font=\scriptsize,anchor=west] at (2.4,-.03) {after $b$};
  \node at (9.25,3.55) {Circulation response};
  \draw[->,black!70,line width=.5pt] (6.9,.5)--(12.15,.5) node[right] {$M$};
  \draw[->,black!70,line width=.5pt] (6.9,.5)--(6.9,3.12) node[above] {$\omega$};
  \draw[black!70,line width=.65pt] (6.9,.85)--(11.85,2.9);
  \draw[popblue,dashed,line width=.4pt] (8.3,.5)--(8.3,1.43)--(6.9,1.43);
  \draw[popred,dashed,line width=.4pt] (10.6,.5)--(10.6,2.38)--(6.9,2.38);
  \draw[popblue,fill=white,line width=.8pt] (8.3,1.43) circle[radius=2.3pt];
  \fill[popred] (10.6,2.38) circle[radius=2.3pt];
  \node[font=\scriptsize,below] at (8.3,.5) {$M_0$};
  \node[font=\scriptsize,below] at (10.6,.5) {$M_s$};
  \node[font=\scriptsize,left] at (6.9,1.43) {$\omega(M_0)$};
  \node[font=\scriptsize,left] at (6.9,2.38) {$\omega(M_s)$};
\end{tikzpicture}
\caption{Inward transport increases the weighted occupancy from $M_0$ to $M_s$, and hence the angular speed when $\omega_1>0$. Rotation preserves $M$, so the two orders differ by the angle $t\omega_1(M_s-M_0)$.}
\label{fig:crowd}
\end{figure}
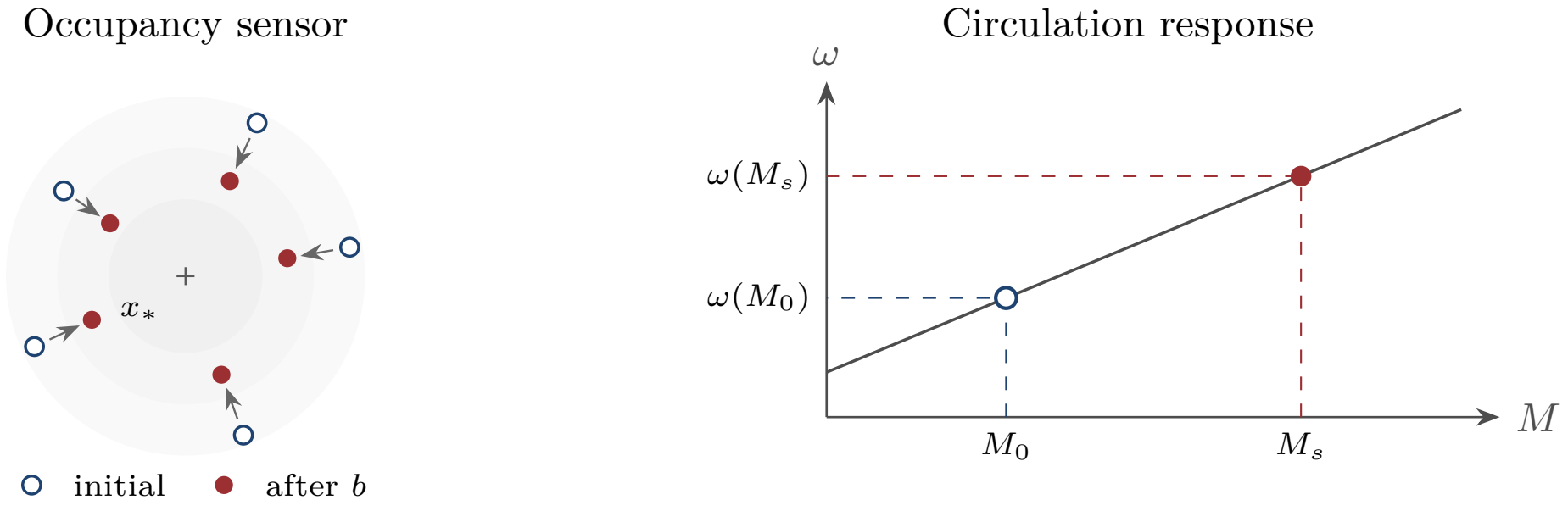

Since $\partial_\mu M(\mu)(y)=\nabla\eta(y)$, we have
$$
\partial_\mu f(x,\mu)(y) =-\kappa_1(x-x_*)\otimes\nabla\eta(y),\qquad \partial_\mu g(x,\mu)(y) =\omega_1J(x-x_*)\otimes\nabla\eta(y).
$$
By radial symmetry, we have $\nabla\eta(y)\cdot J(y-x_*)=0$ and $\nabla\eta(y)\cdot(y-x_*)=-\vert y-x_*\vert^2\eta(y)/\sigma^2$. Consequently,
$$
\br{f}{g}(x,\mu) =\kappa(M(\mu))\omega_1Q(\mu)J(x-x_*),\qquad Q(\mu)=\int\frac{\vert y-x_*\vert^2}{\sigma^2}\eta(y)\,d\mu(y).
$$
Although the frozen spatial bracket vanishes, the Lie-Lions bracket is a rotation about $x_*$, whose angular speed is $\omega_1$ times the rate $\mathscr{L}_fM=\kappa(M)Q$ at which gathering increases the occupancy, whereas circulation preserves the occupancy, i.e., $\mathscr{L}_gM=0$. Gathering before circulating therefore produces a larger rotation of the crowd (see Figure~\ref{fig:crowd}).

If $\omega_1>0$ and $\mu\ne\delta_{x_*}$, then $Q(\mu)>0$. For $\mu=\delta_z$ with $z\ne x_*$, the bracket is nonzero at its atom. Hence the global flows do not commute by Theorem~\ref{thm:commutation}. Writing $M_s=M(\Phi_s^f\mu)$ and $M_0=M(\mu)$, the two compositions have the same final radii and an exact relative angle $\theta(s,t)=t[\omega(M_s)-\omega(M_0)] =t\omega_1(M_s-M_0)$. Indeed, the radial aggregation flow is equivariant under rotations, and the sensor is rotation invariant. Since $\frac{d}{ds}M_s\vert_{s=0}=\mathscr{L}_fM(\mu)=\kappa(M_0)Q(\mu)$, we have $\partial_s\theta(0,t)=t\kappa(M_0)\omega_1Q(\mu)$, which is $t$ times the angular speed of the flow of $\br{f}{g}$ about $x_*$. The two final laws differ whenever the transported law is not invariant under this relative rotation. In particular, a rotation-invariant initial law is unchanged by the angular difference, although the bracket field may be nonzero. Since $f,g\in\SL$, Theorem~\ref{thm:quantitative} gives a cubic remainder for the short-time commutator in $W_2$.

\section{Concluding remarks}\label{sec:conclusion}

We have introduced the Lie-Lions bracket of nonlocal vector fields, $\br{f}{g}=[f,g]_x+[f,g]_\mu$, which adds to the classical Lie bracket in the space variable a correction due to the evolution of the law. Under first-order regularity and linear growth assumptions, its vanishing characterizes the commutation of the flows on $\Ptwo$ (Theorem~\ref{thm:commutation}). The particle realizations preserve the bracket exactly, for every number of agents, and jointly determine the nonlocal field. For the Lie algebras $\SL$ and $\mathcal{B}_L$ of polynomial nonlocal fields, we have established commutator and second-order Baker-Campbell-Hausdorff formulas with cubic remainders in $W_2$, uniformly in the number of agents (Theorem~\ref{thm:quantitative} and Proposition~\ref{prop:bounded-kernels}). The Lie-Lions derivatives of these fields act on polynomial observables through a Leibniz rule, and their action at a fixed law is represented by a tangent vector in the sense of Otto. The examples of Section~\ref{sec:examples} show that the vanishing of the classical bracket at a fixed law is neither necessary nor sufficient for commutation. They also show that the Lie-Lions bracket describes concrete effects: the dependence of a rotation on the inertia or on the occupancy of a population, the phase drift induced by synchronization, the transmission of a turn along a hierarchy, and effective three-body interactions created by alternating pairwise interactions.

\medskip

We conclude with two directions for future research.

\paragraph{Controllability.}
Let us consider the control system $\partial_t\mu+\operatorname{div}\big((u_1f+u_2g)\mu\big)=0$. By Theorem~\ref{thm:quantitative}, the composition $\Phi_{-t}^g\circ\Phi_{-t}^f\circ\Phi_t^g\circ\Phi_t^f$ moves the population along the flow of $\br{f}{g}$ during the time $t^2$, up to a remainder of order three in $W_2$, uniformly in the number of agents. The bracket may provide directions of motion which are not combinations of $f$ and $g$ at the current law. In the flock model of Section~\ref{subsec:flock-turning}, at a law $\mu_{\boldsymbol\theta}$ with $\theta_1=\theta_2$, neither $f$ nor $g$ changes the heading of the follower, whereas $\br{f}{g}$ does as soon as $\sin(\theta_3-\theta_2)\neq0$. Iterated brackets are used in \cite{GoncharovaPogodaevStaritsyn2026} to control finitely many observables. Obtaining controllability results with estimates of the control cost that are uniform in the number of agents, by combining such rank conditions with Theorem~\ref{thm:quantitative}, is an open problem.

\paragraph{Measure differential equations.}
Measure differential equations \cite{Piccoli2019} replace the velocity field by a probability measure on the tangent bundle, which allows concentration and splitting of mass. Extending the Lie-Lions bracket and the estimates of this paper to that setting is a natural continuation.

\appendix

\section{Appendix}
\subsection{Differentiation of integral observables}\label{app:calculus}
We prove the differentiability of the integral observables used in Proposition~\ref{prop:operator}.

\begin{lemma}
\label{lem:iterability}
Let $f$ satisfy Assumption~\ref{ass:regularity}, and let $\varphi\in C^2(\R^d)$ have bounded Hessian. Then $G_\varphi(\mu)=\int\nabla \varphi(x)\cdot f(x,\mu)\,d\mu(x)$ has a continuously Fr\'echet differentiable lift on $\mathcal{H}$. Its Lions derivative is
\begin{equation}\label{iterability-derivative}
\partial_\mu G_\varphi(\mu)(z) = D^2\varphi(z)f(z,\mu)+\partial_xf(z,\mu)^T\nabla \varphi(z) + \int\partial_\mu f(x,\mu)(z)^T\nabla \varphi(x)\,d\mu(x),
\end{equation}
$\mu$-almost everywhere in $z$. The last integral is well defined in $L^2(\mu;\R^d)$, and for every $R>0$ there is $K_{\varphi,R}$ such that $\Vert \partial_\mu G_\varphi(\mu)\Vert_{L^2(\mu)}\leq K_{\varphi,R}$ whenever $M_2(\mu)\leq R$.
\end{lemma}

\begin{proof}
We set $q(x,X)=\nabla \varphi(x)\cdot\widetilde f(x,X)$. This is a $C^1$ map on $\R^d\times \mathcal{H}$. By identifying the dual of $\mathcal{H}$ with $\mathcal{H}$, we consider $\partial_Xq(x,X)$ as an element of $\mathcal{H}$. Since $\vert \nabla \varphi(x)\vert \leq C_\varphi(1+\vert x\vert)$, Assumption~\ref{ass:regularity} implies, for $\Vert X\Vert_2\leq R$,
\begin{equation}\label{integral-derivative-bounds}
\vert \partial_xq(x,X)\vert \leq K_R(1+\vert x\vert), \qquad \Vert \partial_Xq(x,X)\Vert_2\leq K_R(1+\vert x\vert)^2.
\end{equation}
Here and below the constants may also depend on $\varphi$ and $f$. The lifted functional is $\widetilde G_\varphi(X)=\mathbb{E}[q(X,X)]$, where the first argument is evaluated pointwise. It is finite by \eqref{field-growth}.

We define $A(X)\in\mathcal{H}$ by
\begin{equation}\label{integral-candidate}
A(X)=\partial_xq(X,X)+\mathbb{E}[\partial_Xq(X,X)],
\end{equation}
where the expectation of the $\mathcal{H}$-valued map $\omega\mapsto\partial_Xq(X(\omega),X)$ is a Bochner integral in $\mathcal{H}$. By \eqref{integral-derivative-bounds}, the first term belongs to $\mathcal{H}$ and the Bochner integral is well defined. Moreover, $\Vert A(X)\Vert_2$ is bounded on every ball in $\mathcal{H}$.

We first prove that $A$ is norm-continuous. Let $X_k\to X$ in $\mathcal{H}$, with all norms bounded by $R$. Along any subsequence, a further subsequence converges pointwise almost everywhere. Joint continuity then gives $\partial_xq(X_k(\omega),X_k)\rightarrow \partial_xq(X(\omega),X)$ and $\partial_Xq(X_k(\omega),X_k)\rightarrow \partial_Xq(X(\omega),X)$ almost everywhere, with convergence in $\mathcal{H}$ in the second expression. Because $X_k\to X$ in $L^2$, the family $(1+\vert X_k\vert)^2$ is uniformly integrable. The first estimate in \eqref{integral-derivative-bounds} and Vitali's theorem imply $\Vert \partial_xq(X_k,X_k)-\partial_xq(X,X)\Vert_2\rightarrow0$. The second estimate and the same theorem imply $\mathbb{E}\Vert \partial_Xq(X_k,X_k)-\partial_Xq(X,X)\Vert_2\rightarrow0$. Every subsequence has a further subsequence with these properties, so both limits hold for the original sequence. Formula~\eqref{integral-candidate} now gives $\Vert A(X_k)-A(X)\Vert_2\to0$.

For $Y\in \mathcal{H}$, we apply the fundamental theorem of calculus to $\theta\mapsto q(X(\omega)+\theta Y(\omega),X+\theta Y)$ and integrate over $\omega$. Interchanging these integrals is justified by \eqref{integral-derivative-bounds}: the two derivative terms are bounded by constants times $(1+\vert X\vert +\vert Y\vert)\vert Y\vert $ and $(1+\vert X\vert +\vert Y\vert)^2\Vert Y\Vert_2$, respectively. Both are integrable. Consequently
$$
\widetilde G_\varphi(X+Y)-\widetilde G_\varphi(X) =\int_0^1\mathbb{E}\big[A(X+\theta Y)\cdot Y\big]\,d\theta.
$$
The continuity of $A$ implies Fr\'echet differentiability with $D\widetilde G_\varphi(X)Y=\mathbb{E}[A(X)\cdot Y]$, and this derivative is continuous.

It remains to identify $A(X)$. The first term in \eqref{integral-candidate} is
$$
\partial_xq(X,X) =D^2\varphi(X)f(X,\mu)+\partial_xf(X,\mu)^T\nabla \varphi(X), \qquad \mu=\mu_X.
$$
For the second term, \eqref{lions-representation} gives $\partial_Xq(x,X)=\partial_\mu f(x,\mu)(X)^T\nabla\varphi(x)$, so that $\mathbb{E}[\partial_Xq(X,X)]=\eta(X)$, with $\eta(z)=\int\partial_\mu f(x,\mu)(z)^T\nabla \varphi(x)\,d\mu(x)$. This is an $L^2(\mu)$-valued Bochner integral. Indeed,
$$
\Vert \eta\Vert_{L^2(\mu)} \leq\int\vert \nabla \varphi(x)\vert \Vert \partial_\mu f(x,\mu)\Vert_{L^2(\mu)}\,d\mu(x) \leq K_R\int(1+\vert x\vert)^2\,d\mu(x).
$$
The constant includes the norm equivalence stated after Assumption~\ref{ass:regularity}. This proves \eqref{iterability-derivative}. Its stated bound follows from the last estimate and the linear growth of the spatial terms.
\end{proof}

\subsection{Weighted particle estimates}\label{app:quantitative}
We prove the uniform estimates used for Theorem~\ref{thm:quantitative}. All derivative norms below are multilinear operator norms for the specified norm on the particle space.

\begin{lemma}
\label{lem:weighted-particles}
Let $f\in\SL$ and let $r\geq0$ be an integer. There is $C_r>0$ such that, for every $N$ and every initial configuration $\boldsymbol x=(x_1,\ldots,x_N)$, the particle field $f^N$ satisfies
\begin{equation}\label{weighted-particle-bounds}
\Vert D^j f^N(\boldsymbol y) [\boldsymbol h^1,\ldots,\boldsymbol h^j]\Vert_{\boldsymbol x,*} \leq C_r\prod_{\ell=1}^j \Vert \boldsymbol h^\ell\Vert_{\boldsymbol x,*}, \qquad 0\leq j\leq r,
\end{equation}
whenever $\Vert \boldsymbol y-\boldsymbol x\Vert_{\boldsymbol x,*}<1/4$. For $j=0$, the right-hand side is $C_r$. The constant is independent of $N$ and $\boldsymbol x$.
\end{lemma}

\begin{proof}
We set $a_i=\langle x_i\rangle$. Since $z\mapsto\langle z\rangle$ is $1$-Lipschitz,
\begin{equation}\label{weight-comparison}
\frac{3}{4}a_i\leq\langle y_i\rangle\leq\frac{5}{4}a_i
\end{equation}
on the indicated ball. It suffices to treat a kernel representation with $m$ measure variables, for which
$$
f_i^N(\boldsymbol y) =\frac{1}{N^m}\sum_{j_1,\ldots,j_m=1}^N K(y_i;y_{j_1},\ldots,y_{j_m}).
$$
An order-$j$ directional derivative is a sum of at most $(m+1)^j$ types of terms, each averaged over the same $N^m$ indices. If $s_0$ derivatives act on the distinguished slot and $s_\ell$ on its $\ell$th measure slot, with $\sum_{\ell=0}^m s_\ell=j$, then estimate \eqref{kernel-symbol} and \eqref{weight-comparison} bound the kernel derivative by $C a_i^{1-s_0}\prod_{\ell=1}^m a_{j_\ell}^{-s_\ell}$. Each direction acting on a particle coordinate $k$ has magnitude at most $a_k\Vert \boldsymbol h^\nu\Vert_{\boldsymbol x,*}$. Multiplying these inequalities cancels all the negative powers and gives $Ca_i\prod_{\nu=1}^j\Vert \boldsymbol h^\nu\Vert_{\boldsymbol x,*}$ for each term. This remains valid when some indices coincide, since derivatives are assigned to kernel slots before summing over indices. The normalized sum introduces no factor depending on $N$. Dividing the $i$th component by $a_i$ and taking the maximum proves \eqref{weighted-particle-bounds}.
\end{proof}

\begin{lemma}
\label{lem:finite-dimensional-formulas}
Let $r,M>0$, let $E$ be a finite-dimensional normed space and let $x\in E$. Let $u,v$ be $C^3$ fields on the ball $B(x,r)$ such that $\Vert D^ju(y)\Vert +\Vert D^jv(y)\Vert \leq M$ for every $y\in B(x,r)$ and every $j\in\{0,\ldots,3\}$.
We set $w=Dv\,u-Du\,v$ and $z_t=t(u+v)+t^2w/2$. There are $t_0,C>0$, depending only on $r,M$, such that, for $\vert t\vert \leq t_0$, all the following flows are defined along the trajectories in question and
\begin{align}
\label{finite-commutator} \Vert \Phi_{-t}^v\Phi_{-t}^u \Phi_t^v\Phi_t^u(x)-\Phi_{t^2}^w(x)\Vert &\leq C\vert t\vert^3,\\
\label{finite-bch} \Vert \Phi_t^v\Phi_t^u(x)-\Phi_1^{z_t}(x)\Vert &\leq C\vert t\vert^3.
\end{align}
These constants do not depend on $\dim E$ or on the chosen norm.
\end{lemma}

\begin{proof}
By assumption, $w$, $Dw$ and $D^2w$ are bounded by constants depending only on $M$. We choose $t_0\leq1$ small enough that the four successive $u,v$ legs of duration $\vert t\vert $, the $w$ leg of duration $t^2$, and the $z_t$ leg of duration $1$ starting at $x$ stay in $B(x,r/2)$. This is possible because their speeds are bounded respectively by $M$, a constant depending on $M$, and $C_M\vert t\vert $. The same bounds apply to intermediate trajectories needed below by decreasing $t_0$ once.

For a $C^2$ field $V$ with the corresponding bounds, differentiating its trajectory three times, we get
$$
\frac{d^3}{ds^3}\Phi_s^V(y) =D^2V(\Phi_s^V(y)) [V(\Phi_s^V(y)),V(\Phi_s^V(y))] +DV(\Phi_s^V(y))DV(\Phi_s^V(y)) V(\Phi_s^V(y)).
$$
Taylor's formula with its integral remainder therefore yields, uniformly at the intermediate points,
\begin{equation}\label{flow-second-order}
\Phi_{\varepsilon t}^V(y) =y+\varepsilon t\,V(y) +\frac{t^2}{2}DV(y)V(y)+R, \qquad \Vert R\Vert \leq C_M\vert t\vert^3, \quad \varepsilon\in\{-1,1\}.
\end{equation}
If an intermediate point has an expansion $y=x+tA+t^2Q+\mathrm{O}(\vert t\vert^3)$ with bounded $A,Q$, ordinary Taylor expansion gives
$$
V(y)=V(x)+tDV(x)A+\mathrm{O}(t^2),\qquad DV(y)V(y)=DV(x)V(x)+\mathrm{O}(\vert t\vert).
$$
Substitution in \eqref{flow-second-order} proves the recursion
\begin{equation}\label{composition-recursion}
A\mapsto A+\varepsilon V(x),\qquad Q\mapsto Q+\varepsilon DV(x)A +\frac{1}{2}DV(x)V(x)
\end{equation}
for composing one additional leg. All remainders are bounded by constants depending only on $M,r$, since only four legs are used.

For the commutator in \eqref{finite-commutator}, the successive signed fields, in their order of application, are $u,v,-u,-v$. Starting with $A=Q=0$, recursion \eqref{composition-recursion} gives the successive linear coefficients $u$, $u+v$, $v$, $0$. The successive quadratic coefficients, with all fields and derivatives evaluated at $x$, are $\frac{1}{2}Du\,u$, $\frac{1}{2}Du\,u+Dv\,u+\frac{1}{2}Dv\,v$, $Dv\,u-Du\,v+\frac{1}{2}Dv\,v$, $Dv\,u-Du\,v$. Consequently the commutator equals $x+t^2w(x)+\mathrm{O}(\vert t\vert^3)$. Since $w,Dw$ are uniformly bounded, we obtain $\Phi_{t^2}^w(x)=x+t^2w(x)+\mathrm{O}(t^4)$, which proves \eqref{finite-commutator}.

For the BCH estimate, the first two legs give
\begin{equation}\label{bch-composition-expansion}
\Phi_t^v\Phi_t^u(x) =x+t(u+v) +t^2\big(\frac{1}{2}Du\,u+Dv\,u+\frac{1}{2}Dv\,v\big) +\mathrm{O}(\vert t\vert^3).
\end{equation}
On $B(x,r)$ the field $z_t$ and its first two derivatives are bounded by $C_M\vert t\vert $. The third time derivative of its trajectory is therefore bounded by $C_M\vert t\vert^3$. Applying the same integral Taylor formula over $s\in[0,1]$ yields
\begin{align*}
\Phi_1^{z_t}(x) &=x+z_t(x)+\frac{1}{2}Dz_t(x)z_t(x)+\mathrm{O}(\vert t\vert^3)\\
&=x+t(u+v) +\frac{t^2}{2}\big(w+D(u+v)(u+v)\big) +\mathrm{O}(\vert t\vert^3).
\end{align*}
The quadratic coefficient is $\frac{1}{2}Du\,u+Dv\,u+\frac{1}{2}Dv\,v$, because $w=Dv\,u-Du\,v$. Comparison with \eqref{bch-composition-expansion} proves \eqref{finite-bch}.
\end{proof}

\paragraph{Acknowledgements.}
Vincent Boulard and Amaury Hayat were financially supported by the ANR-Tremplin project StarPDE (ANR-24-ERCS-0010).

\paragraph{AI statement.}
The exposition was revised with the assistance of GPT-6.1 Sol and Claude Opus 5.5, which were also used for bibliographic searches and the preparation of the figures. The authors take full responsibility for the content of this article.

\end{document}